\documentclass[10pt,a4paper,twoside,fleqn]{article}
\usepackage{amsfonts,amssymb,amsmath,amsthm,latexsym,mathrsfs}
\usepackage{a4wide,fancyhdr}
\usepackage{graphicx,epic,eepic}
\usepackage{titlesec}
\titleformat{\section}{\large\bfseries}{\thesection}{1em}{}

\numberwithin{equation}{section}

\newtheorem{theorem}{Theorem}[section]
\newtheorem{lemma}{Lemma}[section]

\newtheorem{proposition}{Proposition}[section]

\newtheorem{corollary}{Corollary}[section]

\theoremstyle{definition}

\begin{document}

\begin{center}
{\bf\large The spectral shift function for planar obstacle scattering at low energy}
\end{center}

\begin{center}
I McGillivray\\
School of Mathematics\\
University of Bristol\\
University Walk\\
Bristol BS8 1TW\\
United Kingdom\\
e maiemg@bristol.ac.uk\\
t + 44 (0)117 3311663\\
f + 44 (0)117 9287999
\end{center}

\medskip

\begin{abstract}
Let $H$ signify the free non-negative Laplacian on $\mathbb{R}^2$ and $H_Y$ the non-negative Dirichlet Laplacian on the complement $Y$ of a nonpolar compact subset $K$ in the plane. We derive the low-energy expansion for the Krein spectral shift function (scattering phase) for the obstacle scattering system $\left\{\,H_Y,\,H\,\right\}$ including detailed expressions for the first three coefficients. We use this to investigate the large time behaviour of the expected volume of the pinned Wiener sausage associated to $K$.
\end{abstract}


\medskip


\section{Introduction}

\noindent Given a nonpolar compact subset $K$ in the plane, we consider the exterior domain $Y$ complementary to $K$. Let $\mathscr{H}$ stand for the Hilbert space $L^2(\mathbb{R}^2,\,m)$ where $m$ signifies the Lebesgue measure. The (non-negative) Laplacian acting in $\mathscr{H}$ will be denoted $H$ while $H_Y$ denotes the (non-negative) Dirichlet Laplacian acting in $L^2(Y,\,m)$. The operator $J$ embeds $L^2(Y,\,m)$ into $\mathscr{H}$ through extending by zero on $Y$. Let $\xi(\lambda)$ stand for the Krein spectral shift function (scattering phase) for the pair $\left\{\,H_Y,\,H\right\}$. The following trace formula then holds
\begin{equation}\label{trace_formula_1}
\mathrm{Tr}\,[\,J g(H_Y) J^*-g(H)\,]=\int_0^\infty g^\prime(\lambda)\xi(\lambda)\,d\lambda
\end{equation}
for any function $g:\,\mathbb{R}\rightarrow\mathbb{R}$ of the form
\[
g(\lambda) = \left\{
\begin{array}{lcl}
f(e^{-\lambda}) & \text{ for } & \lambda>0,\\
0               & \text{ for } & \lambda\leq 0,
\end{array}
\right.
\]
where $f:\,\mathbb{R}\rightarrow\mathbb{R}$ is a continuously differentiable function with $f^\prime\in W^{1,2}(\mathbb{R})$. In this paper, we study the asymptotic behaviour of $\xi(\lambda)$ for small $\lambda$. Our main result is this. 

\medskip

\begin{theorem}\label{expansion_for_xi_of_lambda}
Let $K$ be a nonpolar compact subset of $\mathbb{R}^2$. Let $l\in\mathbb{N}$. There exist $\xi_0^k\in\mathbb{R}$ ($-\infty<k\leq -1$) such that
\[
\xi(\lambda)=\sum_{k=-l}^{-1}
\xi_0^k\,(-\log\,\lambda)^k 
+ o((-\log\,\lambda)^{-l})
\]
as $\lambda\downarrow 0$. The first three coefficients are given by 
\begin{itemize}
\item[{\it (i)}]
$\xi_0^{-1}=1$; 
\item[{\it (ii)}]
$\xi_0^{-2}=C(K) - \log\,4 + 2\gamma$; 
\item[{\it (iii)}]
$\xi_0^{-3}=\left(\,C(K) - \log\,4 + 2\gamma\,\right)^2 - \frac{\pi^2}{3}$.
\end{itemize}
\end{theorem}

\medskip

\noindent The quantity $C(K)$ is related to the Robin constant $R(K)$ of $K$ via the relation $C(K)=-4\pi\,R(K)$; $\gamma$ stands for Euler's constant. We remark that the leading order expansion for $\xi(\lambda)$ has been derived in \cite{HassellZelditch1999}. The counterpart to this result in higher dimensions may be found in \cite{McG2000}, \cite{McG2010}.

\medskip

\noindent The study is motivated in part by its relevance to a problem in probability theory. The pinned Wiener sausage $S(t,\omega)$ refers  to the random set swept out by a compact set $K$ in the plane as it is transported along a Brownian loop $\omega:[0,t]\rightarrow\mathbb{R}^2$ by rigid motion. In detail,
\[
S(t,\omega):=\bigcup_{0\leq s\leq t}(\omega(s)+K)\subseteq\mathbb{R}^2;
\]
its area is denoted $|S(t,\omega)|$. Introduce the expected area via
\[
\gamma(t):=\mathbb{E}_{0,0}^{0,t}|S(t,\omega)|
\]
where $\mathbb{P}_{0,0}^{0,t}$ signifies the Brownian bridge measure on loop space associated to the Laplacian $\Delta$. We derive an asymptotic expansion for the above quantity in the large time r\'{e}gime. 

\medskip

\begin{theorem}\label{main_theorem}
Let $K$ be a nonpolar compact subset of $\mathbb{R}^2$. Let $l\in\mathbb{N}$. Then there exist $\gamma_0^k\in\mathbb{R}$ ($\mathbb{Z}\ni k\leq -1$) such that
\[
\gamma(t)=\sum_{k=-l}^{-1}\gamma_0^k\,t\,(\log\,t)^k + o(t\,(\log\,t)^{-l})
\]
as $t\rightarrow\infty$. The first three coefficients are given by
\begin{itemize}
\item[{\it (i)}]
$\gamma_0^{-1}=4\pi$;
\item[{\it (ii)}]
$\gamma_0^{-2}=4\pi\left\{\,C(K) + \gamma  - \log\,4 \right\}$;
\item[{\it (iii)}]
$\gamma_0^{-3}=4\pi\left\{\,(\,C(K)+\gamma - \log\,4)^2-\frac{\pi^2}{6}\right\}$.
\end{itemize}
\end{theorem}

\medskip

\noindent This result was conjectured in \cite{vandenBerg1997}. This latter work derives the first order asymptotic expansion of $\gamma(t)$ for an arbitrary nonpolar compact $K$. It also obtains the third order asymptotic series as above for the particular case in which $K=K_a$ is a closed disc with radius $a>0$. 

\medskip

\noindent It is interesting to compare the behaviour of $\gamma(t)$ with the related functional $\beta(t):=\mathbb{E}_{0}|S(t,\omega)|$. Here $\omega:\,[0,\infty)\rightarrow\mathbb{R}^2$ is a Brownian path in $\mathbb{R}^2$ and $\mathbb{P}_{0}$ stands for the Wiener measure on path space associated to $\Delta$. In the large time r\'{e}gime
\[
\beta(t)=\sum_{k=-l}^{-1}\beta_0^k\,t\,(\log\,t)^k + o(t\,(\log\,t)^{-l})
\]
with explicit expressions for the first three coefficients according to \cite{LeGall1990}. This last expansion extends the work of \cite{Spitzer1964} which detailed the second order expansion. The series for $\gamma(t)$ and $\beta(t)$ agree to leading order. For the lower order terms we have that
\[
\beta_0^{-2}=4\pi\left\{\,C(K) + 1 + \gamma  - \log\,4 \right\}
\hspace{5mm}
\text{ while}
\hspace{5mm}
\beta_0^{-3}=4\pi\left\{\,(\,C(K)+ 1 + \gamma - \log\,4)^2-\frac{\pi^2}{6}\right\}.
\]

\medskip

\noindent The analogous problem for $\gamma$ in higher dimensions has been treated in  \cite{McG2000}, \cite{McG2010}. This problem originated in the calculation of the specific heat of a quantum system of obstacles $K$ at low temperature \cite{UhlenbeckBeth1936}. 

\medskip

\noindent In Section 2 we introduce the trace formula (\ref{trace_formula_1}) and relate the Krein spectral shift function $\xi(\lambda)$ to the scattering matrix $S(\lambda)$ for the system $\left\{H_Y,\,H\right\}$ via the Birman-Kre\v{\i}n formula. We show that this relation holds in particular on an interval of the form $(0,\delta)$, including the case when $K$ does not have a connected complement. In Section 3 we derive a number of prerequisite results in logarithmic potential theory. 

\medskip

\noindent In order to construct the scattering matrix $S(\lambda)$ it is necessary to invert the operator 
\begin{equation}\label{problem_with_invertibility}
I+R^{(-1)}(\mu-\imath\,0)V
\end{equation}
in $B(\mathscr{H}_{-s})$ for $\lambda>0$ in a neighbourhood of $\lambda=0$; here, $\mu$ relates to $\lambda$ via $\mu=(\lambda+1)^{-1}$. To explain terminology briefly, 
\[
V = J R_Y(-1) J^* - R(-1)
\] 
denotes the difference between the Dirichlet and free resolvents; while 
$R^{(-1)}(\cdot)$ signifies the resolvent of $R(-1)$. Also, $\mathscr{H}_{-s}$ refers to a weighted Hilbert space. The operator in (\ref{problem_with_invertibility}) explodes on the complement of a hyperspace in $\mathscr{H}_{-s}$. This complication is absent in higher dimensions; it presents the salient technical difficulty of the paper. This is tackled in Sections 4 and 6.

\medskip

\noindent Section 6 continues with a small energy expansion of the scattering matrix $S(\lambda)$ in a double-series akin to expansions obtained in \cite{Bolle1988}, \cite{Gesztesy1987}. A lattice-point counting lemma in Section 5 plays a role in establishing summability of the double-series. Section 6 culminates in the proof of the expansion given in Theorem \ref{expansion_for_xi_of_lambda}. The detailed derivation of the coefficients is left to Section 7. The application Theorem \ref{main_theorem} is proved in Section 8. The Appendix includes the proofs of several results from \cite{McG2000}.  
\section{The trace formula}

\noindent{\em The free Laplacian.}  Let $\mathscr{H}$ stand for the complex Hilbert space $L^2(\mathbb{R}^2,m)$ based on Lebesgue measure $m$ with inner product $(\cdot,\cdot)$ linear in the first factor. We refer to the non-negative Laplacian $-\Delta$ in $\mathscr{H}$ by $H$. Its resolvent $R(\zeta):=(H-\zeta)^{-1}$ $(\zeta\in\mathbb{C}\setminus[0,\infty))$ has convolution kernel $k(x;\zeta)$ given by
\begin{equation}\label{definition_of_k_x_zeta}
k(x;\zeta):=\frac{\imath}{4}H^{(1)}_0(\zeta^{1/2}|x|)
\end{equation}
where $H_0^{(1)}$ is the first Hankel function of order $0$. The condition $\Im\zeta^{1/2}>0$ specifies the branch of $\zeta^{1/2}$. For the sake of completeness, we recall that
\[
H^{(1)}_0(z)=1+\frac{2}{\pi}\gamma\imath
-\left\{1+\frac{2}{\pi}\imath(\gamma-1)\right\}\frac{z^2/4}{(1!)^2}
+\left\{1+\frac{2}{\pi}\imath(\gamma-1-\frac{1}{2})\right\}\frac{(z^2/4)^2}{(2!)^2}
+\cdots
\]
\begin{equation}\label{expansion_of_Hankel_function}
+\frac{2}{\pi}\gamma\imath\mathrm{Log}(z/2)\left\{ 1 - \frac{z^2/4}{(1!)^2} + \frac{(z^2/4)^2}{(2!)^2}-\cdots\right\}
\hspace{1cm}
(z\in\mathbb{C}\setminus[0,\,\infty))
\end{equation}
as in \cite{AbramowitzStegun1964} 9.1.3, 9.1.12, 9.1.13. The logarithm $\mathrm{Log}$ refers to the principal branch of the logarithm. 

\medskip

\noindent Let us introduce constants
\begin{equation}\label{formula_for_coefficient_a_j}
a_j  :=  
\left\{
\begin{array}{ll}
(1/2\pi)\,\left(\,\log\,2-\gamma\,\right)+\imath / 4, & j=0,\\
\left\{\,\frac{1}{2\pi}\left(\,\log\,2-\gamma-\sum_{k=1}^j\frac{1}{k}\,\right)+\frac{\imath}{4}\,\right\}\frac{(-1)^j}{4^j(j!)^2}, & j\geq 1;
\end{array}
\right.
\end{equation}
\begin{equation}\label{formula_for_coefficient_b_j}
b_j  := 
\left\{
\begin{array}{ll}
-1/2\pi, & j=0,\\
\frac{(-1)^{j+1}}{2\pi}\frac{1}{4^j(j!)^2}, & j\geq 1;
\end{array}
\right.
\end{equation}
\begin{equation}\label{formula_for_coefficient_c_j}
c_j  :=  \frac{1}{4\pi}\frac{(-1)^{j}}{4^j(j!)^2},\hspace{1cm}j\geq 0. 
\end{equation}
Put
\[
k^0_j(x) = \left\{\,a_j + b_j\,\log\,|x|\,\right\}\,|x|^{2j}
\text{ and }
k^1_j(x) = c_j\,|x|^{2j}
\hspace{1cm}(\,x\in\mathbb{R}^2\setminus\{0\}\,).
\]
Then
\begin{equation}\label{expansion_for_resolvent_kernel_k}
k(x;\zeta)=\sum_{j=0}^\infty\sum_{\varepsilon=0}^1\zeta^j\eta^\varepsilon\,k_j^\varepsilon(x)\hspace{1cm}(x\in\mathbb{R}^2\setminus\{0\})
\end{equation}
with
$
\eta:=-2\,\mathrm{Log}\,\zeta^{1/2}.
$

\medskip

\noindent Define $\langle x\rangle:=(1+|x|^2)^{1/2}$ for $x\in\mathbb{R}^2$. The weighted $L^2$-space $\mathscr{H}_s$ ($s\in\mathbb{R}$) is defined by $\mathscr{H}_s:=\{u:\langle\cdot\rangle^s u \in\mathscr{H}\}$. Considered as Banach spaces, the dual space of $\mathscr{H}_s$ is $\mathscr{H}_{-s}$. We write $\langle\cdot,\cdot\rangle$ for the corresponding duality pairing. 

\medskip

\noindent According to \cite{Agmon1995} Theorem 4.1,
\[
R(\lambda\pm\imath\,0):=\lim_{\varepsilon\downarrow 0}R(\lambda\pm\imath\varepsilon)
\]
exists in $B(\mathscr{H}_{s},\mathscr{H}_{-s})$ for any $s>1/2$ and $\lambda>0$, with convergence in the uniform operator topology. Further,

\medskip

\begin{theorem}\label{expansion_of_R_of_zeta}
Let $l\in\mathbb{N}_0$ and $s>2l+1$. Then for $\zeta\in\mathbb{C}\setminus[0,\infty)$,
\[
\left|\left|\,R(\zeta)-\sum_{j=0}^l\sum_{\varepsilon=0}^1\zeta^j\,\eta^\varepsilon K_j^\varepsilon\,\right|\right|_{B(\mathscr{H}_{s},\mathscr{H}_{-s})}=o(|\zeta|^l)
\]
as $\zeta\rightarrow 0$ where $K^\varepsilon_j$ is the operator with convolution kernel $k_j^\varepsilon$.
\end{theorem}

\medskip

\noindent This is proved in the Appendix (see also Proposition 3.7 in \cite{McG2000}). 

\medskip

\noindent The modified resolvent $R^{(-1)}(\zeta)$ is the resolvent of $R(-1)$. It relates to the resolvent of $H$ via
\[
R^{(-1)}((1+\zeta)^{-1})=-(1+\zeta)\left(I+(1+\zeta)R(\zeta)\right),
\hspace{1cm}
\zeta\in\mathbb{C}\setminus[0,\infty).
\]

\noindent The auxiliary Hilbert space $L^2(S^1,\sigma)$ is denoted by $\mathfrak{h}$. Let $U:\mathscr{H}\rightarrow L^2((0,\infty);\mathfrak{h})$ be the spectral representation of $H$. Then for any $u\in\mathscr{H}_s$ and $\lambda>0$,
\[
U(\lambda)u(\omega)=(1/\sqrt{2})\mathscr{F}u(\lambda^{1/2}\omega)\hspace{1cm}(\omega\in S^1)
\] 
provided that $s>1$. Here, $\mathscr{F}$ stands for the Fourier transform
\[
\mathscr{F}u(\xi)=\widehat{u}(\xi):=(2\pi)^{-1}\int_{\mathbb{R}^2}e^{-\imath\,\xi\cdot x}\,u(x)\,m(dx).
\]
The following lemma is proved in the Appendix (see also \cite{McG2000} Lemma 3.9).

\begin{lemma}\label{expansion_of_U_lambda}
Fix $l\in\mathbb{N}_0$ and $s> l + 1$. Then
\[
\left|\left|\,U(\lambda) - \sum_{j=0}^l(\imath\,\lambda^{1/2})^j\,U_j\,\right|\right|_{\mathfrak{S}_2(\mathscr{H}_s,\mathfrak{h})}=o(\lambda^{l/2})
\]
as $\lambda\downarrow 0$. The operator $U_j$ has kernel
\begin{equation}\label{kernel_of_u_j}
u_j(\omega, x) = \frac{1}{\sqrt{2}}(2\pi)^{-1}\frac{(-1)^j}{j!}(\omega\cdot x)^j.
\end{equation}
\end{lemma}

\medskip

\noindent Incidentally, the notation $\mathfrak{S}_2(\mathscr{H}_s,\mathfrak{h})$ refers to the collection of operators from $\mathscr{H}_s$ to $\mathfrak{h}$ of Hilbert-Schmidt type.

\medskip

\noindent Let $U^{(-1)}:\mathscr{H}\rightarrow L^2((0,1);\mathfrak{h})$ be the spectral representation of $R(-1)$. Then $U$ and $U^{(-1)}$ are related via
\begin{equation}\label{U_vis_a_vis_modified_U_relation}
U^{(-1)}(\mu)=(\lambda+1)U(\lambda).
\end{equation}
We use the notation
\[
\mu=(\lambda+1)^{-1}
\]
and this is used routinely in the sequel. Let $l\in\mathbb{N}_0$ and $s>l+1$. From Lemma \ref{expansion_of_U_lambda} we derive
\begin{equation}\label{expansion_for_U_minus_1}
U^{(-1)}(\mu) = \sum_{j=0}^{2l}(\imath\,\lambda^{1/2})^j\,U_j^{(-1)} + o(\lambda^l)
\end{equation}
in $\mathfrak{S}_2(\mathscr{H}_s,\mathfrak{h})$ as $\lambda\downarrow 0$ where
\begin{eqnarray}
U^{(-1)}_0 & = & U_0,\nonumber\\
U^{(-1)}_1 & = & U_1,\\\label{formula_for_U_modified_0}
U^{(-1)}_j & = & U_j - U_{j-2}\text{ for }j\geq 2.\nonumber
\end{eqnarray}

\medskip 

\noindent{\em The spectral shift function.} Let $K$ be a nonpolar compact subset of $\mathbb{R}^2$. Its complement will be denoted by $Y$. Let $H_Y$ refer to the non-negative Dirichlet Laplacian on $L^2(Y,m)$. The semigroup difference
\begin{equation}\label{trace_class_semigroup_difference}
Je^{-H_Y}J^* - e^{-H}\in\mathfrak{S}_1(\mathscr{H})
\end{equation}
is trace class \cite{Stollmann1994}. Let $\xi(\lambda, e^{-H_Y}, e^{-H})$ be the spectral shift function for the pair $\left\{ e^{-H_Y}, e^{-H} \right\}$ (\cite{Yafaev1992} Theorem 8.2.1). Define
\begin{equation}\label{spectral_shift_function}
\xi(\lambda)=\xi(\lambda, H_Y, H):=
\left\{
\begin{array}{ll}
-\xi(e^{-\lambda}, e^{-H_Y}, e^{-H}), & \lambda\geq 0,\\
0,                                    & \lambda<0.
\end{array}
\right.
\end{equation}
By \cite{Yafaev1992} Theorem 8.2.1,
\begin{equation}\label{integrability_of_xi}
\xi\in L^1(\mathbb{R};\,e^{-|\lambda|}\,d\lambda).
\end{equation}
By \cite{Yafaev1992} Theorem 8.3.3 and the paragraph following it, we may write
\begin{equation}\label{trace_formula}
\mathrm{Tr}\,[\,Jg(H_Y)J^*-g(H)\,]=\int_0^\infty g^\prime(\lambda)\xi(\lambda)\,d\lambda
\end{equation}
for any function $g:\,\mathbb{R}\rightarrow\mathbb{R}$ of the form
\[
g(\lambda) = \left\{
\begin{array}{lcl}
f(e^{-\lambda}) & \text{ for } & \lambda>0,\\
0               & \text{ for } & \lambda\leq 0,
\end{array}
\right.
\]
where $f:\,\mathbb{R}\rightarrow\mathbb{R}$ is a continuously differentiable function with $f^\prime\in W^{1,2}(\mathbb{R})$. In particular, given $t>3/2$ we can find a continuously differentiable function $f$ such that $f(\lambda)=\lambda^t$ for $0\leq\lambda\leq 1$ and $f^\prime\in W^{1,2}(\mathbb{R})$. We then have
\begin{equation}\label{trace_formula_for_heat_semigroup}
\mathrm{Tr}\,[\,Je^{-tH_Y}J^*-e^{-tH}\,]=-\int_0^\infty t e^{-t\lambda}\xi(\lambda)\,d\lambda
\end{equation}
for $t>3/2$. 
\medskip

\noindent Let $Y_{e}$ stand for the unbounded connected component of $Y$ and set $Y_b:=Y\setminus Y_e$. In case $Y_b\neq\emptyset$, we differentiate between $H_e$ resp. $H_b$, the non-negative Dirichlet Laplacians on $Y_e$ resp. $Y_b$. The spectrum $\sigma(H_b)$ of $H_b$ is discrete. Let $\xi_e(\lambda)=\xi(\lambda,H_{Y_e},H)$ be the spectral shift function for the pair $\left\{\,H_{Y_e},\,H\right\}$. Denote by
\[
N_b(\lambda):=\sum_{\sigma(H_b)\ni\nu<\lambda}m(\nu)
\]
the spectral counting function for $H_b$; here, $m(\nu)$ stands for the geometric multiplicity of $\nu\in\sigma(H_b)$. 

\medskip

\begin{lemma}
It holds that
\begin{itemize}
\item[{\it (i)}]
$
\xi(\lambda) = \xi_e(\lambda) + N_b(\lambda)
$
for a.e. $\lambda>0$;
\item[{\it (ii)}]
$\xi$ admits an a.e.-version that is real analytic on $(0,\infty)\setminus\sigma(H_b)$.
\end{itemize}
Finally, with $\xi$ denoting this version,
\begin{itemize}
\item[{\it (iii)}] $\xi(\lambda)\rightarrow 0$ as $\lambda\downarrow 0$. 
\end{itemize}
\end{lemma}

\noindent{\em Proof.}
First note that both $\xi$ and $\xi_e$ satisfy (\ref{integrability_of_xi}). For $t>3/2$, 
\[
\int_0^\infty e^{-t\lambda}\xi(\lambda)\,d\lambda
= t^{-1}\,\mathrm{Tr}\,[\,e^{-t\,H} - J_e e^{-tH_{Y_e}} J_e^*\,] - t^{-1}\,\mathrm{Tr}\,[\,e^{-t\,H_{Y_b}}\,]
= \int_0^\infty e^{-t\lambda}\left\{\,\xi_e(\lambda) + N_b(\lambda)\,\right\}\,d\lambda
\]
where the Weyl asymptotics of $N_b(\cdot)$ ensure the absolute integrability of the second integrand. Item {\it (i)} follows by the inversion formula for the Laplace-Lebesgue integral (\cite{Widder1946} Theorem VII.6a). Parts {\it (ii)} and {\it (iii)} follow from \cite{JensenKato1978} Lemmas 3.2 and 3.4 and {\it (i)}.  
\qed

\medskip

\noindent{\em The scattering matrix.} In virtue of (\ref{trace_class_semigroup_difference}) the scattering operator $S(e^{-H_Y},e^{-H},J)$ for the pair $\left\{\,e^{-H_Y},e^{-H}\,\right\}$ exists and is unitary on $\mathscr{H}$ by \cite{Yafaev1992} Theorem 6.2.1 and Corollary 2.4.2. By the invariance principle (\cite{Yafaev1992} Theorem 6.2.5), the scattering operator $S(R_Y(-1),R(-1),J)$ exists and $S(R_Y(-1),R(-1),J)=S(e^{-H_Y},e^{-H},J)$. As the scattering operators commute with the corresponding spectral projectors (\cite{Yafaev1992} Theorem 2.1.4 and 1.5.1) we have the representation
\[
S(\lambda,e^{-H_Y},e^{-H},J)=S(\varphi(\lambda),R_Y(-1),R(-1),J)\hspace{1cm}
\text{ a.e. }\lambda>0
\]
where $\varphi:\,(0,1)\rightarrow(0,1);\,\lambda\mapsto(-\log\lambda+1)^{-1}$. By the Birman-Kre\v{\i}n formula (\cite{Yafaev1992} Theorem 8.4.1),
\begin{eqnarray}\label{Birman_Krein_formula}
e^{2\pi\imath\,\xi(\lambda)} & = & \mathrm{Det}(S(e^{-\lambda},e^{-H_Y},e^{-H},J))\nonumber\\
& = & \mathrm{Det}(S(\mu,R_Y(-1),R(-1),J))\hspace{1cm}\text{ a.e. }\lambda>0.
\end{eqnarray}
We now derive a representation formula for $S(\mu,R_Y(-1),R(-1),J)$.

\medskip

\noindent Set
\[
V:=JR_Y(-1)J^* - R(-1).
\]
Then

\medskip

\begin{theorem}\label{extension_and_compactness_of_V}
For each $s>0$, $V$ admits a bounded extension from $\mathscr{H}_{-s}$ to $\mathscr{H}_{s}$ and $V=V^*\in\mathfrak{S}_\infty(\mathscr{H}_{-s},\mathscr{H}_s)$ is compact.  
\end{theorem}

\medskip

\noindent This result is proved in the Appendix (see also \cite{McG2000} Theorem 2.1). Given $s>1/2$ and $\mu\in(0,1)$, define
\[
\mathscr{H}_s^\mu:=\left\{\,f=(\,R(-1)-\mu\,)\,u:\,u\in\mathscr{H}_s\,\right\}.
\]

\begin{lemma}\label{properties_of_R(-1)}
Let $s>1/2$ and $\mu\in(0,1)$. We have
\begin{itemize}
\item[{\it (i)}]
$R(-1)\in B(\mathscr{H}_s,\mathscr{H}_s)$;
\item[{\it (ii)}] $\mathscr{H}_s^\mu$ is a proper subspace of $\mathscr{H}_s$;
\item[{\it (iii)}] the identity
\[
R^{(-1)}(\mu\pm\imath\,0)\,(\,R(-1)-\mu\,) = I
\]
holds on $\mathscr{H}_s$.
\end{itemize}
\end{lemma}

\noindent{\em Proof.}
{\it (i)} From the identity
$
\langle x\rangle^s\leq 2^s\,\left\{\,\langle y\rangle^s + \langle x-y \rangle^s\,\right\}
$
obtain
\[
\langle x\rangle^s k(x-y;-1)\langle y\rangle^{-s}\leq k(x-y;-1)+\langle x-y \rangle^sk(x-y;-1)
\hspace{1cm}(x\neq y).
\]
The latter kernel is integrable by \cite{Taylor1996a} 3.6, so defines a bounded convolution operator on $\mathscr{H}$ by Young's inequality \cite{Davies1989}.

\medskip

\noindent{\it (ii)} Let us introduce the Sobolev space  $W_s:=\left\{u:\,\widehat{u}\in\mathscr{H}_s\,\right\}$. Let $\tau:\,W_s\rightarrow L^2(S^1_\lambda,\sigma)$ stand for the restriction mapping (\cite{ReedSimon1975} Theorem IX.39). Then $\tau\,\widehat{f}=0$ for any $f\in\mathscr{H}_s^\mu$. The function $f=e^{-|\cdot|^2/2}\in\mathscr{H}_s$ does not satisfy this property as $\widehat{f}=f$ (\cite{LiebLoss1997} Theorem 5.2).

\medskip

\noindent{\it (iii)} Let $u\in\mathscr{H}_s$ and $f:=(\,R(-1)-\mu\,)u$. Then
\begin{eqnarray*}
\left|\left|\,u - R^{(-1)}(\mu-\imath\,0)f\,\right|\right|_{\mathscr{H}_{-s}}
& = &
\lim_{\varepsilon\downarrow 0}\left|\left|\,u - R^{(-1)}(\mu-\imath\,\varepsilon)f\,\right|\right|_{\mathscr{H}_{-s}}\\
& = & \lim_{\varepsilon\downarrow 0}\left|\left|\,u - R^{(-1)}(\mu-\imath\,\varepsilon)\left[\,R(-1)-(\mu-\imath\,\varepsilon)-\imath\,\varepsilon\,\right]u\,\right|\right|_{\mathscr{H}_{-s}}\\
& = & \lim_{\varepsilon\downarrow 0}\varepsilon\,\left|\left|\, R^{(-1)}(\mu-\imath\,\varepsilon)u\,\right|\right|_{\mathscr{H}_{-s}}\\
& = & 0
\end{eqnarray*}
and similar with the opposite sign. 
\qed

\begin{theorem}\label{absence_of_embedded_eigenvalues}
Let $s>1/2$.
\begin{itemize}
\item[{\it (i)}]
Assume that $Y_b=\emptyset$. For any $\mu\in(0,1)$, the compact operator $V R^{(-1)}(\mu\pm\imath\,0)$ acting in $B(\mathscr{H}_s)$ does not have eigenvalue $-1$.  
\item[{\it (ii)}]
Assume that $Y_b\neq\emptyset$ and $\lambda\not\in\sigma(H_b)$. Then $V R^{(-1)}(\mu\pm\imath\,0)$ acting in $B(\mathscr{H}_s)$ does not have eigenvalue $-1$.  
\item[{\it (iii)}] 
Assume that $Y_b\neq\emptyset$ and $\lambda\in\sigma(H_b)$. Then $V R^{(-1)}(\mu\pm\imath\,0)$ acting in $B(\mathscr{H}_s)$ has eigenvalue $-1$.  
\end{itemize}
\end{theorem}

\noindent{\em Proof.}
{\it (i)} Suppose that $V R^{(-1)}(\mu-\imath\,0)f=-f$ for some $f\in\mathscr{H}_s$. Set $u:=R^{(-1)}(\mu-\imath\,0)f$. Argue as in \cite{DemuthMcGillivray1999} Lemma 4.4 to conclude that $u\in\mathscr{H}$ and that $R_Y(-1)u = \mu\,u$. Put $f:=\Delta\,u + \lambda\,u\in\mathscr{D}^\prime(\mathbb{R}^2)$ with $\lambda:=-1+1/\mu$. Then $f=0$ on $Y$ because $u$ is a weak solution of $\Delta\,u+\lambda\,u=0$ there. By elliptic regularity \cite{Taylor1996a} Proposition 3.9.1, $u$ is smooth on $Y$. Adapting the argument in \cite{Taylor1996b} Lemma 1.2 to the $d=2$ case, conclude that $u$ vanishes on the complement of some ball $B(0,r)$. The unique continuation property (\cite{Hormander1983} Theorem 5.1, for example) ensures that $u$ vanishes throughout $Y$. The proof of {\it (ii)} is similar. 

\medskip

\noindent{\it (iii)} Let $\varphi\in L^2(Y_b)$ be an eigenfunction of $H_b$ corresponding to $\lambda$. Let $u\in\mathscr{H}_s$ be the extension of $\varphi$ by $0$. Then $f:=\left[\,R(-1)-\mu\,\right]u\in\mathscr{H}_s$ by Lemma \ref{properties_of_R(-1)} {\it (i)}. Also, $R^{(-1)}(\mu-\imath\,0)f=u$ by Lemma \ref{properties_of_R(-1)} {\it (iii)}. In an obvious notation,
\[
Vu =\left[\,R_{Y_b}(-1)\oplus R_{Y_e}(-1) - R(-1)\,\right]u
= -\left[\,R(-1)-\mu\,\right]u = -f;
\]
that is, $VR^{(-1)}(\mu-\imath\,0)f=-f$. 
\qed

\medskip

\noindent For $\lambda\in(0,\infty)\setminus\sigma(H_b)$, 
\[
\exists\,(\,I+VR^{(-1)}(\mu+\imath 0))^{-1}\in B(\mathscr{H}_s)
\]
by the Fredholm alternative.
As in \cite{Yafaev1992} Theorem 5.7.1' (with $\mathfrak{G}=\mathscr{H}_s$ for $s>1/2$ and $G:\,\mathscr{H}\rightarrow\mathfrak{G};\,f\mapsto\langle\cdot\rangle^{-2s} f$) the scattering matrix for $\{\,R_Y(-1), R(-1)\,\}$ can be represented
\[
S(\mu,R_Y(-1),R(-1),J)
=
I - 2\pi\imath\,U^{(-1)}(\mu)(\,I+VR^{(-1)}(\mu+\imath 0))^{-1}VU^{(-1)}(\mu)^*,
\hspace{1cm}\text{ a.e. }\lambda\in(0,\infty)\setminus\sigma(H_b).
\]

\medskip

\noindent Let $S(\cdot)$ stand for the (adjoint) scattering matrix for $\{H,\,H_Y\}$,
\[
S(\lambda)=I + 2\pi\imath\,U^{(-1)}(\mu)V(I+R^{(-1)}(\mu-\imath 0)V)^{-1}U^{(-1)}(\mu)^*
\in B(\mathfrak{h})
\]
with $\lambda\in(0,\infty)\setminus\sigma(H_b)$. The $t$-matrix is characterised by the relation $S(\lambda)=I+T(\lambda)$. From (\ref{Birman_Krein_formula}),
\begin{equation}\label{Birman_Krein_formula}
e^{-2\pi\imath\xi(\lambda)}=\mathrm{Det}\,S(\lambda)\hspace{1cm}\text{ a.e. }\lambda>0.
\end{equation}
\section{Some logarithmic potential theory}

\noindent For brevity, we use the notation $g^\lambda$ to stand for the resolvent operator $R(-\lambda)$ with $\lambda>0$; $g^\lambda(\cdot)$ stands for the coresponding convolution kernel. If $\lambda=0$ the notation $g$ is sometimes used. For $z\in\mathbb{C}\setminus(-\infty,0]$ define
\begin{equation}\label{definition_of_b}
b(z):=H^{(1)}_0(z) - 1 - \frac{2\imath}{\pi}\left\{\,\mathrm{Log}(z/2)+\gamma\,\right\}.
\end{equation}
Given $0<\delta<1$ there exists a finite constant $c$ such that
\begin{equation}\label{estimate_for_b}
\left|\,b(z)\,\right|\leq C\,|z|^2(-\log\,|z|)\text{ for }|z|\leq\delta.
\end{equation}
The logarithmic potential kernel is defined by 
\begin{equation}\label{definition_of_k}
k(x):= (1/2\pi)\,\log(\,1/|x|\,)
\hspace{1cm}(x\in\mathbb{R}^2\setminus\{0\}).
\end{equation}
From (\ref{definition_of_k_x_zeta}), the kernel $g^\lambda(\cdot)$ may be decomposed in terms of (\ref{definition_of_k}) and (\ref{definition_of_b}) as 
\begin{equation}\label{decomposition_of_g_lambda}
g^\lambda(x)=a_0-\imath\pi c_0 - c_0\log\lambda + k(x) + r(\lambda^{1/2}|x|)
\hspace{1cm}(x\in\mathbb{R}^2\setminus\{0\})
\end{equation}
with $a_0$ and $c_0$ as in (\ref{formula_for_coefficient_a_j}) and (\ref{formula_for_coefficient_c_j}), and
\[
r(x):=(\imath/4)b(\imath x).
\]
Fix a unit vector $u$ in $\mathbb{R}^2$. The regularised resolvent kernel $k^\lambda(\cdot)$ is given by 
\begin{equation}\label{formula_for_k_lambda}
k^\lambda(x) := g^\lambda(x)-g^\lambda(u)\hspace{1cm}(x\in\mathbb{R}^2\setminus\{0\}).
\end{equation}
From (\ref{decomposition_of_g_lambda}),
\[
g^\lambda(u)=a_0 - \imath\pi c_0 - c_0\,\log\,\lambda + r(\lambda^{1/2}),
\]
and hence
\begin{equation}\label{formula_for_k_lambda_minus_k}
k^\lambda(x)-k(x)=r(\lambda^{1/2}|x|) - r(\lambda^{1/2})\hspace{1cm}(x\in\mathbb{R}^2\setminus\{0\}).
\end{equation}
The operators with convolution kernels $k^\lambda(\cdot)$ resp. $k(\cdot)$ will be denoted by $k^\lambda$ resp. $k$. We use the notation $r^\lambda$ to refer to the operator with convolution kernel $r(\lambda^{1/2}|x|)$.

\medskip

\begin{corollary}\label{convergence_of_k_lambda_to_k_in_B_H_s_H_minus_s}
Let $s>1$. Then $k^\lambda\rightarrow k$ in $B(\mathscr{H}_s,\mathscr{H}_{-s})$ as $\lambda\downarrow 0$. 
\end{corollary}

\noindent{\em Proof.}
We may write
\[
k^\lambda - k = g^\lambda - \eta\,K_0^1 - K_0^0 - r(\lambda^{1/2})\langle\cdot,1\rangle1
\]
with $\eta$ as before given by $\eta=-2\,\mathrm{Log}(-\lambda)^{1/2}=-\log\,\lambda-\imath\,\pi$. 
Now apply Theorem \ref{expansion_of_R_of_zeta} and (\ref{estimate_for_b}).
\qed

\begin{theorem}\label{convergence_of_r_lambda_locally_uniformly}
Let $s>1$. Then
\begin{itemize}
\item[{\it (i)}] $k:\,\mathscr{H}_s\rightarrow C(\mathbb{R}^2)$;
\item[{\it (ii)}] $r^\lambda:\,\mathscr{H}_s\rightarrow C(\mathbb{R}^2)$ for each $\lambda>0$;
\item[{\it (iii)}]
$g^\lambda(u)\,r^\lambda f \rightarrow 0$ locally uniformly on $\mathbb{R}^2$ as $\lambda\downarrow 0$ for each $f\in\mathscr{H}_{s}$.
\end{itemize}
\end{theorem}

\noindent{\em Proof.}
{\it (i)} Define
\[
k_1 f(x) := -(1/2\pi)\int_{B(x,1)}\log\,|x-y|\,f(y)\,m(dy),
\]
and likewise for $k_2$ but with $B(x,1)$ replaced by its complement $B(x,1)^c$. Note that
\[
\left|\,\int_{B(x,r)}f(y)\,m(dy)\,\right|
\leq
\sqrt{\pi}\,\left|\left|\,f\,\right|\right|_{\mathscr{H}_s}\,r
\]
for $f\in\mathscr{H}_s$. As in \cite{Adams1996} Lemma 3.1.1 (b),
\[
k_1 f(x) = (1/2\pi)\int_0^1\int_{B(x,r)}f(y)\,m(dy)\,\frac{dr}{r}
+ (1/2\pi)\int_{B(x,1)}f(y)\,m(dy).
\]
A dominated convergence argument shows that $k_1f$ is continuous on $\mathbb{R}^2$. As for $k_2$, we have
\[
\chi_{B(x,1)^c}(y)\log\,|x-y| \leq \tau + |y|\hspace{1cm}(y\in\mathbb{R}^2)
\] 
for all $x\in B(0,\tau)$ $(\tau>0)$. Continuity of $k_2 f$ follows by another appeal to dominated convergence. 

\medskip

\noindent{\it (ii)} Write
\[
r^\lambda = g^\lambda  - k + \left\{\,r(\lambda^{1/2})- g^\lambda(u)\,\right\}\langle\cdot,1\rangle\,1.
\]
Now $\mathscr{H}_s\subseteq L^p(\mathbb{R}^2)$ for each $p>1$. By the Sobolev embedding \cite{Adams1996} Theorem 1.2.4, $g^\lambda:\,\mathscr{H}_s\rightarrow C(\mathbb{R}^2)$. This and {\it (i)} establish the claim. 

\medskip

\noindent{\em (iii)} Pick $0<\alpha<1/3$. For each $x\in\mathbb{R}^2$ introduce sets
\[
\begin{array}{lll}
A^\lambda_1 & := & \{\,y\in\mathbb{R}^2:\,|y-x|\leq 1\},\\
A^\lambda_2 & := & \{\,y\in\mathbb{R}^2:\,1<|y-x|\leq\lambda^{-\alpha}\},\\
A^\lambda_3 & := & \{\,y\in\mathbb{R}^2:\,\lambda^{-\alpha}<|y-x|\leq\delta\,\lambda^{-1/2}\},\\
A^\lambda_4 & := & \{\,y\in\mathbb{R}^2:\,|y-x|>\delta\,\lambda^{-1/2}\},
\end{array}
\]
for $\lambda$ sufficiently small (where the $x$-dependence has been suppressed for the sake of legibility). Define
\[
r^\lambda_j f(x) := \int_{A^\lambda_j}r(\lambda^{1/2}|x-y|) f(y)\,m(dy)
\hspace{0.5cm}(\,x\in\mathbb{R}^2\,)\hspace{0.5cm}(j=1,2,3,4).
\]
Fix $\tau>0$. For $0<\lambda<\delta^2$,
\[
\left|\,r^\lambda_1 f(x)\,\right|
\leq
(C/4)\,\lambda
\left\{\,
\int_{|x-y|\leq 1}
\left(\,\log\lambda^{1/2}\,|x-y|\,\right)^2\,m(dy)\,
\right\}^{1/2}\,
\left|\left|\,f\,\right|\right|_{\mathscr{H}_s}
\]
by (\ref{estimate_for_b}). In particular, $g^\lambda(u)\,r^\lambda_1 f\rightarrow 0$ uniformly on $B(0,\tau)$ as $\lambda\downarrow 0$.

\medskip

\noindent Again from (\ref{estimate_for_b}), for $\lambda>0$ small,
\[
\left|\,r^\lambda_2 f(x)\,\right|
\leq
(C/4)\,\lambda
\left\{\,
\int_{1<|x-y|\leq\lambda^{-\alpha}}\,|x-y|^4
\left[\,
(1/2)(\log\lambda)^2+2\,(\log\,|x-y|\,)^2\right]\,m(dy)
\right\}^{1/2}\,
\left|\left|\,f\,\right|\right|_{\mathscr{H}_s}.
\]
Choose $0<\eta<1/\alpha-3$. An estimate of the form $\left|\,\log\,|x-y|\,\right|\leq c_\eta\,|x-y|^\eta$ holds on $A^\lambda_2$. Also,
\[
\left\{\,\int_{1<|x-y|\leq\lambda^{-\alpha}}|x-y|^{4+2\,\eta}\,m(dy)\,\right\}^{1/2}
\leq\left\{\frac{\pi}{3+\eta}\right\}^{1/2}\lambda^{-\alpha(3+\eta)}.
\]
So $g^\lambda(u)\,r^\lambda_2 f\rightarrow 0$ uniformly on $B(0,\tau)$ as $\lambda\downarrow 0$.

\medskip

\noindent The kernel $r(\lambda^{1/2}|x-y|)$ is bounded by a constant $c^\prime$ (say) on $A^\lambda_3$. Thus, for any $x\in B(0,\tau)$,
\begin{eqnarray*}
\left|\,r^\lambda_3 f(x)\,\right|
& \leq &
c^\prime\,\int_{\lambda^{-\alpha}<|x-y|\leq\delta\,\lambda^{-1/2}}\,|f(y)|\,m(dy)\\
& \leq &
c^\prime\,\int_{\lambda^{-\alpha}<|x-y|}\,|f(y)|\,m(dy)\\
& \leq & c^\prime\,\int_{B(0,\lambda^{-\alpha}-\tau)^c}\,|f(y)|\,m(dy)\\
& \leq & c^\prime\,\left|\left|\,f\,\right|\right|_{\mathscr{H}_s}\,
\left\{\,\int_{B(0,\lambda^{-\alpha}-\tau)^c}\langle y\rangle^{-2s}\,m(dy)\,\right\}^{1/2}
\end{eqnarray*}
provided $\lambda$ is sufficiently small. The weight function $\langle\cdot\rangle^{-2s}$ is integrable because $s>1$. Thus $g^\lambda(u)\,r^\lambda_3 f\rightarrow 0$ uniformly on $B(0,\tau)$ as $\lambda\downarrow 0$.

\medskip

\noindent By (\ref{decomposition_of_g_lambda}),
\[
r^\lambda_4 f(x) = \int_{|x-y|>\delta\lambda^{-1/2}}
\left\{\,g^\lambda(x-y)-a_0+\imath\pi c_0 + c_0\,\log\lambda - k(x-y)\,\right\}\,f(y)\,m(dy).
\]
An estimate of the form (\ref{estimate_for_H_1_0_of_z}) gives
\begin{eqnarray*}
\left|\,\int_{|x-y|>\delta\lambda^{-1/2}}
g^\lambda(x-y)\,f(y)\,m(dy)\,\right| & \leq & (c/4)\,\lambda^{-1/4}\,\int_{|x-y|>\delta\lambda^{-1/2}}|x-y|^{-1/2}\,|f(y)|\,m(dy)\\
 & \leq &
(c/4)\,\delta^{-1/2}\,\lambda^{1/2}\int_{|x-y|>\delta\lambda^{-1/2}}\,|f(y)|\,m(dy).
\end{eqnarray*}
Uniform convergence can be derived in a way similar to the treatment of $r^\lambda_3$. 

\medskip

\noindent Choose $\eta>0$ such that $s-\eta>1$. For $x\in B(0,\tau)$ and $\lambda$ small,
\begin{eqnarray*}
\left|\,\int_{|x-y|>\delta\lambda^{-1/2}}
k(x-y)\,f(y)\,m(dy)\,\right| & \leq & c_\eta\,\int_{|x-y|>\delta\lambda^{-1/2}}|x-y|^{\eta}\,|f(y)|\,m(dy)\\
 & \leq &
c_\eta^\prime\,\int_{|x-y|>\delta\lambda^{-1/2}}\langle y\rangle^{\eta}\,|f(y)|\,m(dy)\\
 & \leq &
c_\eta^\prime\,\int_{B(0,\delta\lambda^{-1/2}-\tau)^c}\langle y\rangle^{\eta}\,|f(y)|\,m(dy)\\
& \leq & c_\eta^\prime\,\left|\left|\,f\,\right|\right|_{\mathscr{H}_s}\,
\left\{\,\int_{B(0,\delta\lambda^{-1/2}-\tau)^c}
\langle y\rangle^{-2(s-\eta)}\,m(dy)\,\right\}^{1/2}
\end{eqnarray*}
for appropriate constants $c_\eta$, $c_\eta^\prime$. The remaining terms in $r^\lambda_4$ can be dealt with using similar analysis. Consequently, $g^\lambda(u)\,r^\lambda_4 f\rightarrow 0$ uniformly on $B(0,\tau)$ as $\lambda\downarrow 0$.
\qed

\medskip

\noindent Let $M=(\Omega,\mathscr{M},X_t,\mathbb{P}_{x})$ be Brownian motion on $\mathbb{R}^2$ with transition function $p(t,\,\cdot)$ given by
\[
p(t,x)=(4\pi t)^{-1}e^{-|x|^2/4t}\hspace{1cm}(t>0)
\]
(see \cite{PortStone1978} for example). Put $\sigma_K:=\inf\{t>0:\,X_t\in K\,\}$, the first hitting time of $K$. The hitting operators $h^\lambda_K$ are defined by
\[
h^\lambda_K f :=\mathbb{E}_\cdot\left[ e^{-\lambda\,\sigma_K}f(X_{\sigma_K}):\,\sigma_K<\infty\right]
\]
for measurable $f\geq 0$. If $\lambda=0$, we write $h_K$ instead of $h^0_K$. The $\lambda$-potential of $K$ is the function
\[
p^\lambda_K:=h^\lambda_K 1.
\]
In case $\lambda=0$, $p_K:=p_K^0\equiv 1$ by recurrence (\cite{PortStone1978} Proposition 2.9, for example). Set
\begin{equation}\label{definition_of_W_lambda}
W^\lambda_K:=g^\lambda(u)\,\left\{\,1-p^\lambda_K\,\right\}.
\end{equation}
In case $K$ is nonpolar, the limit
\begin{equation}\label{definition_of_W}
W_K:=\lim_{\lambda\downarrow 0}W^\lambda_K
\end{equation}
exists finitely according to \cite{PortStone1978} Theorem 3.4.2. The equilibrium measure $\mu_K$ is the unique probability measure concentrated on $K^r$ whose potential $k\mu_K$ is constant on $K^r$. Here, $K^r$ denotes the set of regular points for $K$ (\cite{PortStone1978} 2.3). This constant value is the Robin constant $R(K)$ of $K$. We use the notation
\[
C(K):=-4\pi R(K).
\]
From \cite{PortStone1978} Theorem 3.4.12,
\begin{equation}\label{identity_linking_R_and_W}
k\,\mu_K = R(K) - W_K.
\end{equation}
In particular, $W_K\in\mathscr{H}_{-s}$ for any $s>1$. The $1$-capacity of $K$ is denoted $C_1(K)$ (see \cite{Fukushima1980}). It holds that
\begin{equation}\label{formula_for_C_1_of_K}
C_1(K)=\langle\, 1, p^1_K\,\rangle.
\end{equation}

\medskip

\noindent We use the notation $g^\lambda_K$ ($\lambda>0$) to refer to the $\lambda$-potential operator with kernel
\[
g^\lambda_K(x,y):=\int_0^\infty e^{-\lambda\,t}q_K(t,x,y)\,dt\hspace{1cm}(x,y\in\mathbb{R}^2);
\]
to clarify,
$
q_K(t,x,y)=p(t,x-y) - r_K(t,x,y)
$
with
\[
r_K(t,x,y):=\mathbb{E}_x[\,p(t-\sigma_K,X(\sigma_K)-y):\,\sigma_K<t\,]
\]
as in \cite{PortStone1978} 2.5. We have that $g^\lambda_K f=R_Y(-\lambda)f$ $m$-a.e. on $Y$ for $f\in L^2(Y)$. If $\lambda=0$ the notation $g_K$ is sometimes used.

\medskip

\noindent As in \cite{PortStone1978} 3.4, the fundamental identities for logarithmic potentials read
\begin{equation}\label{fundamental_identity_for_logarithmic_potentials_lambda_positive}
k^\lambda = g^\lambda_K + h^\lambda_K\,k^\lambda-\langle\,\cdot\,,\,1\,\rangle\,W^\lambda_K,
\end{equation}
\begin{equation}\label{fundamental_identity_for_logarithmic_potentials_lambda_zero}
k         = g_K + h_K\,k-\langle\,\cdot\,,\,1\,\rangle\,W_K.
\end{equation}

\medskip

\begin{lemma}\label{resolvent_applied_to_potential}
For any $\lambda,\,\mu \geq 0$ with $\lambda+\mu>0$ we have that
\[
g_K^\lambda p_K^\mu=-(\mu-\lambda)^{-1}\left\{\,p_K^\mu - p_K^\lambda\,\right\}.
\]
\end{lemma}

\noindent{\em Proof.}
The result follows as in the proof of \cite{McG2000} Proposition 4.13.
\qed

\medskip

\noindent Let $T_r:=\inf\{t>0:\,X_t\not\in B(0,r)\}$ stand for the first exit time of  $B(0,r)$. The notation $B(0,r)$ signifies the open ball with centre at the origin and radius $r>0$.

\medskip

\begin{lemma}\label{expected_exit_time_from_ball}
For each $r>0$,
\[
\mathbb{E}_0\left[e^{-T_r}\right]=1/I_0(r)
\]
where $I_0$ stands for the modified Bessel function of order zero. Moreover,
\[
\mathbb{E}_0\left[e^{-T_r}\right]\sim(2\pi\,r)^{1/2}\,e^{-r}
\]
as $r\rightarrow\infty$. 
\end{lemma}

\noindent{\em Proof.}
A direct computation leads to the identity. The asymptotic behaviour follows as in \cite{Taylor1996a} 3.6.
\qed

\begin{proposition}\label{matrix_element_W_with_p_1}
We have that
\[
\langle\,W_K,p_K^1\,\rangle = 1.
\]
\end{proposition}

\noindent{\em Proof.}
In a similar way to the proof of \cite{PortStone1978} Proposition 3.4.4, we write
\begin{equation}\label{identity_from_Port_and_Stone}
\langle\,k(\cdot,y)-k(0,y)\,1, p_K^1\,\rangle
=
g_K p_K^1(y) + \langle\,h_K[k(\cdot,y)-k(0,y)\,1], p_K^1\,\rangle
-
\langle\,W_K,p_K^1\,\rangle
\hspace{1cm}(y\in\mathbb{R}^2)
\end{equation}
The function 
\[
k(x,y)- k(0,y)=-\frac{1}{2\pi}\log\frac{|x-y|}{|y|}\hspace{1cm}(x\in\mathbb{R}^2\setminus\{y\})
\]
converges to zero uniformly on compacts as $y\rightarrow\infty$. From Lemma \ref{resolvent_applied_to_potential}, 
\[
g_K p_K^1(y)=1-p_K^1(y).
\]
This latter converges to $1$ in the limit $y\rightarrow\infty$ by Lemma \ref{expected_exit_time_from_ball}. Therefore, the expression on the right-hand side of (\ref{identity_from_Port_and_Stone}) converges to $1-\langle\,W_K,p_K^1\,\rangle$.

\medskip

\noindent We now show that the left-hand side vanishes in the limit. Let $\varepsilon>0$ and choose $0<\delta<1$ with the property that
\[
|\log(1+\tau)| < 2\pi\,\varepsilon\text{ for any }\tau\in\mathbb{R}\text{ with }|\tau|<\delta.
\]
For any $x,y\in\mathbb{R}^2$ with $|x|<\delta\,|y|$ we have 
\begin{equation}\label{estimate_on_k}
|\,k(x,y) - k(0,y)\,|<\varepsilon.
\end{equation}
The left-hand side in (\ref{identity_from_Port_and_Stone}) may be written, 
\[
\langle\,k(\cdot,y)-k(0,y)\,1, p_K^1\,\rangle
=
\int_{\mathbb{R}^2}\left\{\,k(x,y)-k(0,y)\,\right\} p_K^1(x)\,m(dx).
\]
Decompose $\mathbb{R}^2$ into the disjoint union $\mathbb{R}^2=F_y\,\dot{\cup}\,G_y$ with
\[
F_y:=\{\,x\in\mathbb{R}^2:\,|x| < \delta\,|y|\,\}
\text{ and }
G_y:=\{\,x\in\mathbb{R}^2:\,|x| \geq \delta\,|y|\,\}.
\]
The integral over $F_y$ is bounded by $\varepsilon\,\text{Cap}_1(K)$ in modulus by (\ref{estimate_on_k}). Re-write the integral over $G_y$ as
\[
-k(0,y)\int_{G_y}p_K^1(x)\,m(dx)
+\int_{A_y}k(x,y)p_K^1(x)\,m(dx)
+\int_{B_y}k(x,y)p_K^1(x)\,m(dx)
\]
where
\[
A_y:=\{\,x\in\mathbb{R}^2:\,|x|\geq\delta\,|y|,\,|x-y|\leq 1\,\}
\text{ and }
B_y:=\{\,x\in\mathbb{R}^2:\,|x|\geq\delta\,|y|,\,|x-y|>1\,\}.
\]
The first two integrals vanish in the limit by Lemma \ref{expected_exit_time_from_ball}; for the last, use in addition the estimate
\[
|\log\,|x-y|\,|\leq |x|+|y|\text{ for }|x-y|>1.
\]
\qed

\begin{lemma}\label{convergence_of_W_lambda_and_R_Y_minus_lambda}
Let $s>1$. Then
\begin{itemize}
\item[{\it (i)}]
$W_K^\lambda\rightarrow W_K$ in $\mathscr{H}_{-s}$ as $\lambda\downarrow 0$;
\item[{\it (ii)}]
$R_Y(-\lambda)\rightarrow R_Y(0)$ strongly in $B(\mathscr{H}_{s},\mathscr{H}_{-s})$ as $\lambda\downarrow 0$.
\end{itemize}
\end{lemma}

\noindent{\em Proof.}
From (\ref{fundamental_identity_for_logarithmic_potentials_lambda_positive}), (\ref{fundamental_identity_for_logarithmic_potentials_lambda_zero})  and Lemma \ref{resolvent_applied_to_potential} we have
\[
\begin{array}{lll}
k^\lambda p^1_K & = & \frac{1}{1-\lambda}\{\,p^\lambda_K - p^1_K\,\}
+ h^\lambda_K k^\lambda p^1_K - C_1(K) W^\lambda_K,\\
k p^1_K & = & 1-p^1_K 
+ h_K k p^1_K - C_1(K) W_K.
\end{array}
\]
Consequently,
\[
[k^\lambda - k]p^1_K = \frac{1}{1-\lambda}\{\,p^\lambda_K-1\,\}+\frac{\lambda}{1-\lambda}
-\frac{\lambda}{1-\lambda}p^1_K + h_K^\lambda[k^\lambda-k]p^1_K + (h^\lambda_K-h_K)kp^1_K
-C_1(K)\{W^\lambda_K-W_K\}.
\]
{\it (i)} follows with the help of Corollary \ref{convergence_of_k_lambda_to_k_in_B_H_s_H_minus_s} and the dominated convergence theorem. 

\medskip

\noindent From (\ref{fundamental_identity_for_logarithmic_potentials_lambda_positive}) resp. (\ref{fundamental_identity_for_logarithmic_potentials_lambda_zero}) it can be seen that $g^\lambda_K$ resp. $g_K$ map $\mathscr{H}_s$ boundedly into $\mathscr{H}_{-s}$ for any $s>1$. Further,
\[
g^\lambda_K - g_K = k^\lambda - k + h^\lambda_K(\,k^\lambda - k\,)
+ (\,h^\lambda_K - h_K\,) k -\langle\cdot,\,1\rangle\,\left\{\, W^\lambda_K - W_K\,\right\}.
\]
The claim in {\it (ii)} now follows from Corollary \ref{convergence_of_k_lambda_to_k_in_B_H_s_H_minus_s}, Theorem \ref{convergence_of_r_lambda_locally_uniformly} and {\it (i)} above. We use the relation (\ref{formula_for_k_lambda_minus_k}).    
\qed

\begin{proposition}\label{convergence_of_matrix_element_to_1}
It holds that
\[
\lim_{\lambda\downarrow 0}
g^\lambda(u)\left\{\,1-\langle\,W^\lambda_K,p_K^1\,\rangle\,\right\}= R(K).
\]
\end{proposition}

\noindent{\em Proof.}
Applying (\ref{fundamental_identity_for_logarithmic_potentials_lambda_positive}) to the equilibrium measure $\mu_K$ we obtain the identity
\[
k^\lambda\mu_K = h^\lambda_K k^\lambda\mu_K - W^\lambda_K
\]
with the help of \cite{PortStone1978} Theorem 4.4.3 as $\mu_K$ is a probability measure with support in $K^r$. We derive
\[
\langle\,W^\lambda_K,p_K^1\,\rangle =
\langle\,h^\lambda_Kk^\lambda\mu_K,p_K^1\,\rangle - \langle\,k^\lambda\mu_K,p_K^1\,\rangle.
\]
In virtue of (\ref{identity_linking_R_and_W}) we have
\[
\langle\,W_K,p_K^1\,\rangle =
R(K)\,\langle\,1,p_K^1\,\rangle - \langle\,k\mu_K,p_K^1\,\rangle.
\]
Using Proposition \ref{matrix_element_W_with_p_1} we proceed,
\begin{eqnarray*}
1-\langle\,W^\lambda_K,p_K^1\,\rangle
& = &  
R(K)\,\langle\,1,p_K^1\,\rangle
- \langle\,h^\lambda_K[k^\lambda-k]\mu_K,p_K^1\,\rangle
-\langle\,h^\lambda_Kk\mu_K,p_K^1\,\rangle
+ \langle\,[k^\lambda-k]\mu_K,p_K^1\,\rangle\\
 & = & 
\frac{R(K)}{g^\lambda(u)}\langle\,W^\lambda_K,p_K^1\,\rangle
- \langle\,h^\lambda_K[k^\lambda-k]\mu_K,p_K^1\,\rangle
+ \langle\,[k^\lambda-k]\mu_K,p_K^1\,\rangle.
\end{eqnarray*}
The result now follows from Lemma \ref{convergence_of_W_lambda_and_R_Y_minus_lambda} {\it (i)}, (\ref{formula_for_k_lambda_minus_k}), (\ref{estimate_for_b}), Theorem \ref{convergence_of_r_lambda_locally_uniformly} {\it (iii)} and duality.
\qed

\section{Construction of inverse operators}

\noindent Let $X$ be a complex Banach space with dual space $X'$ and duality pairing $\langle\,\cdot\,,\,\cdot\,\rangle$. Let $B(X)$ stand for the collection of bounded linear operators on $X$. The notation $A^\times$ stands for the adjoint operator of $A\in B(X)$. 

\begin{lemma}\label{inverse_of_rank_one_perturbation}
Assume that $A\in B(X)$ is bijective with inverse $B$. Let $y\in X$, $f\in X'$ and $\sigma\in\mathbb{C}$. Define $A_\sigma$ to be the rank-one perturbation of $A$ given by 
\[
A_\sigma:=A+\sigma\langle\,\cdot,f\,\rangle\,y.
\]
Suppose that
\begin{equation}\label{definition_of_alpha}
\alpha:=1+\sigma\,\langle\,By,f\,\rangle\neq 0.
\end{equation}
Then $A_\sigma$ is bijective and has inverse given by
\[
B_\sigma = B - \alpha^{-1}\sigma\,\langle\,\cdot,B^\times f\,\rangle\,By.
\]
\end{lemma}

\noindent{\em Proof.}
Verify by direct computation that $A_\sigma\,B_\sigma = B_\sigma\,A_\sigma = I$. 
\qed

\begin{lemma}\label{abstract_convergence_of_inverses_result}
Let $\delta>0$. Suppose that $(A_\lambda)_{\lambda\in(0,\delta)}$ is a family of operators in $B(X)$. Assume that
\begin{itemize}
\item[{\it (i)}]
$A_\lambda\rightarrow A$ strongly as $\lambda\downarrow 0$ for some $A\in B(X)$;
\item[{\it (ii)}] for each $\lambda\in(0,\delta)$, $A_\lambda$ is bijective  with inverse $B_\lambda$;
\item[{\it (iii)}]
$B_\lambda\rightarrow B$ strongly as $\lambda\downarrow 0$ for some $B\in B(X)$.
\end{itemize}
Then $A$ is bijective and has inverse $B$. 
\end{lemma}

\noindent{\em Proof.}
Given $x\in X$ write, for example,
\[
B\,A\,x-x = (\,B-B_\lambda\,)\,A\,x+B_\lambda\,(\,A - A_\lambda\,)\,x.
\]
By the uniform boundedness principle there exists a finite constant $c$ such that $||B_\lambda||\leq c<\infty$ for all $\lambda\in(0,\delta)$. Take limits on the right-hand side using {\it (i)} and {\it (iii)} to see that $BAx-x=0$.  
\qed

\begin{lemma}\label{simplification_of_I_minus_RV}
For any $\lambda>0$,
\begin{itemize}
\item[{\it (i)}] $\left[\,I-R^{(-1)}_Y((1-\lambda)^{-1})V\,\right]1
=g^\lambda(u)^{-1}W_K^\lambda+\lambda\,p_K^\lambda$;
\item[{\it (ii)}] $\left[\,I-VR^{(-1)}_Y((1-\lambda)^{-1})\,\right]p_K^1
=-V\left[\,g^\lambda(u)^{-1}W_K^\lambda+\lambda\,p_K^\lambda\,\right]$.
\end{itemize}
\end{lemma}

\noindent{\em Proof.} 
Item {\it (ii)} follows from {\it (i)} via the identity $p_K^1 = -V1$. This last  follows from Lemma \ref{resolvent_applied_to_potential}. Again by this lemma,
\begin{eqnarray*}
\left[\,I-R^{(-1)}_Y((1-\lambda)^{-1})V\,\right]1 &=& 1 + R^{(-1)}_Y((1-\lambda)^{-1})p_K^1\\
 &=& 1 - (1-\lambda)\left\{\,I+(1-\lambda)\,R_Y(-\lambda)\,\right\}p_K^1\\
 &=& 1 - (1-\lambda)\,p_K^1 - (1-\lambda)^2\,g^\lambda_K\,p_K^1\\
 &=& 1 - (1-\lambda)\,p_K^1+(1-\lambda)\left\{\,p_K^1 - p_K^\lambda\,\right\}\\
 &=& g^\lambda(u)^{-1}W_K^\lambda+\lambda\,p_K^\lambda.
\end{eqnarray*}
\qed

\noindent Define
\begin{equation}\label{definition_of_A}
\begin{array}{llll}
A         & := & I-\left[\,I+k\,\right]\,V, & \\
A_\lambda & := & I+R^{(-1)}((1-\lambda)^{-1})V-g^\lambda(u)\,\langle\,\cdot\,,\,p_K^1\,\rangle\,1 & (\lambda>0).
\end{array}
\end{equation}

\begin{proposition}\label{properties_of_A_lambda}
Let $s>1$. Then 
\begin{itemize}
\item[{\it (i)}]
$A\in B(\mathscr{H}_{-s})$; 
\item[{\it (ii)}]
$A_\lambda\in B(\mathscr{H}_{-s})$ for any $\lambda>0$; 
\item[{\it (iii)}]
$A_\lambda\rightarrow A$ strongly in $B(\mathscr{H}_{-s})$ as $\lambda\downarrow 0$;
\item[{\it (iv)}]
$
A\,W_K=R(K)\,1
$.
\item[{\it (v)}]
$A^\times V W_K = -R(K)p^1_K$.
\end{itemize}
\end{proposition}

\noindent{\em Proof.}
Statement {\it (i)} flows from \cite{McG2000} Theorem 2.1 and Lemma 3.1. From Theorem \ref{expansion_of_R_of_zeta} and \cite{McG2000} Lemma 3.1, $R(-\lambda)\in B(\mathscr{H}_{s},\mathscr{H}_{-s})$; {\it (ii)} now follows. For {\it (iii)}, we may write
\[
A_\lambda=I-(1-\lambda)\,\left[\,I+(1-\lambda)\,k^\lambda\,\right] V
+ \lambda\,(\lambda - 2)\,g^\lambda(u)\langle\,\cdot\,,\,p^1_K\,\rangle\,1
\]
with the help of (\ref{formula_for_k_lambda}). Thus,
\[
A_\lambda - A = \lambda\,V + \left[\,k - k^\lambda\,\right]\,V
+
\lambda\,(2-\lambda)\,\left\{\,k^\lambda V - g^\lambda(u)\langle\,\cdot\,,\,p^1_K\,\rangle\,1\,\right\}.
\]
The strong convergence follows from Corollary \ref{convergence_of_k_lambda_to_k_in_B_H_s_H_minus_s}. 

\medskip

\noindent As for the identity {\it (iv)}, from Lemma \ref{resolvent_applied_to_potential} and the first resolvent identity, we derive 
\[
\begin{array}{lll}
V p^\lambda_K & = & (1-\lambda)^{-1}\{\,p^\lambda_K - p^1_K\,\} - R(-1)p^\lambda_K,\\
R(-\lambda)Vp^\lambda_K & = & (1-\lambda)^{-1}\left\{\,R(-1)p^\lambda_K - R(-\lambda)p^1_K\,\right\};
\end{array}
\]
the second flowing from the first. With their help, a computation leads to the identity
\[
A_\lambda W_K^\lambda = 
g^\lambda(u)\left\{1 - \langle\,W^\lambda_K,\,p^1_K\,\rangle\,\right\}\,1
- \lambda\,g^\lambda(u)\,p^1_K
+ \lambda(\lambda - 1)g^\lambda(u)R(-\lambda)p^1_K.
\]
By (\ref{formula_for_k_lambda}) and (\ref{formula_for_C_1_of_K}),
\[
\lambda(\lambda - 1)g^\lambda(u)\,R(-\lambda)p^1_K 
= 
\lambda(\lambda - 1)g^\lambda(u)\,
\left\{
k^\lambda p^1_K + g^\lambda(u)\,C_1(K)\,1\,\right\} \rightarrow 0\text{ in }\mathscr{H}_{-s}\text{ as }\lambda\downarrow 0
\]
by Corollary \ref{convergence_of_k_lambda_to_k_in_B_H_s_H_minus_s}. This shows that
\begin{equation}\label{convergence_of_W_lambda_A_lambda}
A_\lambda W^\lambda_K \rightarrow R(K)\,1\text{ in }\mathscr{H}_{-s}\text{ as }\lambda\downarrow 0
\end{equation}
by Proposition \ref{convergence_of_matrix_element_to_1}. Finally, write
\[
A\,W_K - R(K)\,1
=
(\,A-A_\lambda\,)\,W_K + A_\lambda\,(\,W_K - W^\lambda_K\,)
+ A_\lambda W^\lambda_K - R(K)\,1
\]
and use (\ref{convergence_of_W_lambda_A_lambda}), Lemma \ref{convergence_of_W_lambda_and_R_Y_minus_lambda} {\it (i)}, and {\it (iii)}. Of course, the family $(A_\lambda)_{\lambda\in(0,1)}$ is bounded by the uniform boundedness principle.

\medskip

\noindent Lastly,
$
A^\times VW_K = V A W_K = -R(K)p^1_K
$
by {\it (iv)}.
\qed

\begin{lemma}\label{inverse_of_A_lambda_result}
Let $s>1$. Assume that $R(K)\neq 0$. Then there exists $\delta>0$ such that for each $\lambda\in(0,\delta)$ the operator $A_\lambda$ is bijective with inverse $B_\lambda$ given by
\begin{equation}\label{inverse_of_A_lambda}
B_\lambda=I-R^{(-1)}_Y((1-\lambda)^{-1})V
-\frac{1}{g^\lambda(u)\alpha_\lambda}
\langle\,\cdot,V[W_K^\lambda+\lambda g^\lambda(u) p_K^\lambda]\,\rangle\,\left(\,W_K^\lambda+\lambda g^\lambda(u) p_K^\lambda\right)
\end{equation}
where
\[
\alpha_\lambda =
1-\langle W^\lambda_K, p^1_K\rangle + \lambda g^\lambda(u)\langle p^\lambda_K, p^1_K\rangle.
\]
\end{lemma}

\noindent{\em Proof.}
The counterpart $\alpha_\lambda$ of (\ref{definition_of_alpha}) reads
\[
\alpha_\lambda:=1-g^\lambda(u)\langle\left[\,I-R^{(-1)}_Y((1-\lambda)^{-1})V\,\right]1,\,p^1_K\,\rangle
=1-\langle W^\lambda_K, p^1_K\rangle + \lambda g^\lambda(u)\langle p^\lambda_K, p^1_K\rangle
\]
after simplification using Lemma \ref{simplification_of_I_minus_RV}.
By Proposition \ref{convergence_of_matrix_element_to_1}, 
$
g^\lambda(u)\alpha_\lambda\rightarrow R(K)$ as $\lambda\downarrow 0
$.
Consequently, there exists $\delta>0$ such that $\alpha_\lambda\neq 0$ for $\lambda\in(0,\delta)$. By Lemma \ref{inverse_of_rank_one_perturbation}, $A_\lambda$ is bijective with inverse as in (\ref{inverse_of_A_lambda}); again, after making use of Lemma \ref{simplification_of_I_minus_RV}. 
\qed

\noindent Define
\[
B_0 := I + [\,I+R_Y(0)\,]\,V.
\]

\begin{lemma}\label{properties_of_B_0}
The following identities hold:
\begin{itemize}
\item[{\it (i)}]
$B_01 = 0$;
\item[{\it (ii)}]
$B_0^\times p^1_K = 0$.
\end{itemize}
\end{lemma}

\noindent{\em Proof.}
By Lemma \ref{resolvent_applied_to_potential},
\[
B_0 1 = 1 + [\,I+R_Y(0)\,] V1
      = 1 - [\,I+R_Y(0)\,]p^1_K
      = 1 - p^1_K - \left\{\, 1 - p^1_K\,\right\}
      = 0
\]
giving {\it (i)}. For {\it (ii)},
$
B_0^\times p^1_K = - B_0^\times V 1 = - V B_0 1 = 0
$.
\qed

\noindent In case $R(K)\neq 0$, define
\begin{equation}\label{definition_of_B}
B := B_0 - R(K)^{-1}\langle\,\cdot,\,VW_K\,\rangle\,W_K.
\end{equation}

\begin{proposition}\label{B_lambda_converges_to_B}
Let $s>1$.  Assume that $R(K)\neq 0$. Then 
\begin{itemize}
\item[{\it (i)}] 
$B_\lambda\rightarrow B$ strongly in $B(\mathscr{H}_{-s})$ as $\lambda\downarrow 0$;
\item[{\it (ii)}]
$A$ is bijective in $B(\mathscr{H}_{-s})$ with inverse $B$ as in (\ref{definition_of_B});
\item[{\it (iii)}]
$B1 = R(K)^{-1}\,W_K$;
\item[{\it (iv)}]
$B^\times=I + V\,\left[\,I+R_Y(0)\,\right] - R(K)^{-1}\,\langle\cdot,\,W_K\,\rangle\,VW_K$;
\item[{\it (v)}]
$B^\times p_K^1 = -R(K)^{-1}\,VW_K$.
\end{itemize}
\end{proposition}

\noindent{\em Proof.}
{\it (i)} follows from Proposition \ref{convergence_of_matrix_element_to_1} and Lemma \ref{convergence_of_W_lambda_and_R_Y_minus_lambda}. This together with Lemmas \ref{abstract_convergence_of_inverses_result} and \ref{inverse_of_A_lambda_result}, and Proposition \ref{properties_of_A_lambda} {\it (iii)} yield {\it (ii)}. To see {\it (iii)} use Proposition \ref{properties_of_A_lambda} {\it (iv)}. For {\it (v)} use the identity $B^\times\,V=V\,B$ and {\it (iii)}. 
\qed

\begin{lemma}\label{B_0_A_on_M}
Let $s>1$. Set
\begin{eqnarray*}
\mathscr{M} & := & \left\{\,u\in\mathscr{H}_{-s}:\,\langle u, p^1_K\rangle = 0\,\right\},\\
\mathscr{W} & := & \left\{\,u\in\mathscr{H}_{-s}:\,\langle u, VW_K\rangle = 0\,\right\}.
\end{eqnarray*}
Then
\begin{itemize}
\item[{\it (i)}]
$B_0 A = I$ on $\mathscr{M}$;
\item[{\it (ii)}]
$A B_0 = I$ on $\mathscr{W}$.
\end{itemize}
\end{lemma}

\noindent{\em Proof.}
Note that
$
A_\lambda= I + R^{(-1)}((1-\lambda)^{-1})V
$
on $\mathscr{M}$ for each $\lambda>0$.
Set
\[
B_{0,\lambda}:=I - R^{(-1)}_Y((1-\lambda)^{-1})V
\]
on $\mathscr{H}_{-s}$. Then $B_{0,\lambda}\rightarrow B_0$ strongly in $B(\mathscr{H}_{-s})$ as $\lambda\downarrow 0$. By the second resolvent identity and density of $\mathscr{H}\cap\mathscr{M}$ in $\mathscr{M}$, $B_{0,\lambda}A_\lambda = I$ on $\mathscr{M}$. By Proposition \ref{properties_of_A_lambda} {\it (iii)} and as in Lemma \ref{abstract_convergence_of_inverses_result} we obtain $B_0A=I$ on $\mathscr{M}$. This establishes {\it (i)}. 

\medskip

\noindent From Lemma \ref{simplification_of_I_minus_RV} {\it (ii)}, 
$
B_{0,\lambda}^\times p^1_K = - g^\lambda(u)^{-1}VW_K^\lambda - \lambda V p^\lambda_K
$.
Hence,
\[
g^\lambda(u)\langle B_{0,\lambda}u, p^1_K\rangle = \langle u, -VW_K^\lambda-\lambda g^\lambda(u)Vp^\lambda_K\rangle
\rightarrow-\langle u, VW_K\rangle = 0
\text{ as }\lambda\downarrow 0
\]
for $u\in\mathscr{W}$. Therefore,
$
A_\lambda B_{0,\lambda}u = u 
- 
g^\lambda(u)\langle B_{0,\lambda}u, p^1_K\rangle
\rightarrow u
$
as $\lambda\downarrow 0$. Now use strong convergence to obtain {\it (ii)}.
\qed
\section{A lattice-point counting lemma}

\noindent We require a simple lattice-point counting lemma. Let us make the following definitions. For $n\in\mathbb{N}$ and $\mathbb{Z}\ni k<0$ set
\begin{eqnarray*}
A(n,k) & := & \left\{\,x\in\mathbb{Z}^n:\,x_j\leq 1\text{ for }j=1,2,\ldots,n\text{ and }
\sum_{j=1}^nx_j = k\,\right\};\\
a(n,k) & := & \mathrm{Card}(A(n,k)).
\end{eqnarray*}

\begin{lemma}\label{lattice_point_counting_lemma}
For $n\in\mathbb{N}$ and $\mathbb{Z}\ni k<0$, it holds that
\[
a(n,k) \leq a(n)\left\{\,|k|+(3/2)n\,\right\}^{n-1}.
\]
The constant $a(n)$ is given by
\[
a(n) = \frac{\sqrt{n}}{\alpha(n-1)}
\frac{2^{n-1}}{(n-1)!}.
\]
\end{lemma}

\medskip

\noindent  Here, $\alpha(n)$ stands for the volume of the unit ball $B(0,1)$ in $\mathbb{R}^n$; it is understood that $\alpha(0)=1$. 

\medskip

\noindent{\em Proof.}
First, notice that $a(1,k)=1$. For $n=2,3,\ldots$ and $r\leq(3/2)n$ set
\begin{eqnarray*}
H(n,r) & := & \left\{\,x\in\mathbb{R}^n:\,x_j\leq 3/2\text{ for }j=1,2,\ldots,n\text{ and }
\sum_{j=1}^nx_j = r\,\right\},\\
h(n,r) & := & \sigma(H(n,r)),
\end{eqnarray*}
where $\sigma$ stands for surface area measure. We claim that
\begin{equation}\label{formula_for_h_n_r}
h(n,r) = \frac{\sqrt{n}}{(n-1)!}\left\{\,|r|+(3/2)n\,\right\}^{n-1}
\end{equation}
for $r<0$. To see this, introduce the set
\[
S(n,r) := \left\{\,x\in\mathbb{R}^n:\,x_j\geq 0\text{ for }j=1,2,\ldots,n\text{ and }
\sum_{j=1}^n x_j = r\,\right\}
\]
for $r\geq 0$ and $n\geq 2$. As in \cite{Thorpe1979} (for example) its surface area is given by
\[
s(n,r) := \sigma(S(n,r)) = \sqrt{n}\,r^{n-1}/(n-1)!.
\]
Since
\[
H(n,r)=(3/2)\,(1,\ldots,1) - S(n,-r+(3/2)n),
\]
the formula (\ref{formula_for_h_n_r}) follows. 

\medskip

\noindent To prove the lemma, note the inclusion
\[
\dot{\bigcup}_{x\in A(n,k)}B(x,1/2)\cap H(n,k)\subseteq H(n,k)
\]
where the left-hand side is a disjoint union. Computing surface area using (\ref{formula_for_h_n_r}) yields the claim.
\qed

\section{Asymptotics of the spectral shift function}

\noindent The result below follows from Theorem \ref{expansion_of_R_of_zeta}; the method of proof is similar to that used in the proof of \cite{McG2000} Lemma 4.7. 

\medskip

\begin{lemma}\label{expansion_of_R_minus_one_zeta}
Let $l\in\mathbb{N}_0$ and $s>2l+1$. Then for $\zeta\in\mathbb{C}\setminus[0,\infty)$,
\[
\left|\left|\,R^{(-1)}((1+\zeta)^{-1})-\sum_{j=0}^l\sum_{k=0}^1\zeta^j\,\eta^k A_j^k\,\right|\right|_{B(\mathscr{H}_{s},\mathscr{H}_{-s})}=o(|\zeta|^l)
\]
as $\zeta\rightarrow 0$. The coefficients are given by
\begin{eqnarray*}
A_0^1 &=& -K_0^1,\\
A_0^0 &=& -I-K_0^0,\\
A_1^1 &=& -2\,K_0^1 - K_1^1,\\
A_1^0 &=& -I - 2\,K_0^0 - K_1^0,\\
A_j^k &=& -K_j^k - 2\,K_{j-1}^k - K_{j-2}^k\,
\text{ for }j\geq 2\text{ and } k\in\left\{0,\,1\right\}.
\end{eqnarray*}
\end{lemma}

\medskip

\noindent In the context of the last Lemma, we may write
\begin{equation}\label{expansion_for_I_plus_R_V}
I+R^{(-1)}((1+\zeta)^{-1})V = I - \left[\,I+k\,\right]\,V+\sigma\,\langle\,\cdot\,,\,p_K^1\,\rangle\,1
+
\sum_{j=1}^l\sum_{k=0}^1\zeta^j\eta^k A_j^k\,V + o(|\zeta|^l)
\end{equation}
in $B(\mathscr{H}_{-s})$ as $\mathbb{C}\setminus[0,\infty)\ni\zeta\rightarrow 0$. We use the shorthand
\[
\sigma = a + b\,\eta
\]
where $a = a_0$ and $b=c_0$. Define
\[
A_\sigma:=A + \sigma\,\langle\,\cdot\,,\,p_K^1\,\rangle\,1
\]
with $A$ as in (\ref{definition_of_A}).

\begin{proposition}\label{expression_for_B_sigma}
Let $s>1$. Then for small $\zeta\in\mathbb{C}\setminus[0,\infty)$, $A_\sigma\in B(\mathscr{H}_{-s})$ is bijective with inverse given by
\begin{equation}\label{inverse_of_B_sigma}
B_\sigma:= B_0 + \sum_{k=-\infty}^{-1}\theta_k\eta^{k}\langle\,\cdot\,,\,VW_K\,\rangle\,W_K
\end{equation}
where
\begin{equation}\label{formula_for_theta_k}
\theta_k := (-1)^{-k}( 1 / b )\,\left(\,\frac{R(K)+a}{b}\,\right)^{-(k+1)}
\text{ for }k=-1,-2,\ldots.
\end{equation}
Moreover, there exists a finite constant $c$ such that
\[
\left|\left|\,B_\sigma\,\right|\right|_{B(\mathscr{H}_{-s})}
\leq c < \infty
\]
for small $\zeta\in\mathbb{C}\setminus[0,\infty)$.
\end{proposition}

\medskip

\noindent For later use, we introduce the quantity 
\[
\theta:= \left|\frac{R(K)+a}{b}\right|.
\]
Note that $\theta$ is invertible; in fact, $\theta\geq\pi$ for all values of $R(K)\in\mathbb{R}$.

\medskip

\noindent{\em Proof.}
We first treat the case $R:=R(K)\neq 0$. By Lemma \ref{B_lambda_converges_to_B} {\it (iii)} and Proposition \ref{matrix_element_W_with_p_1},
\[
\alpha_\sigma:=1+\sigma\,\langle\,B1,\,p^1_K\,\rangle
= 1 + \sigma / R.
\]
The above quantity is non-zero for small $\zeta\in\mathbb{C}\setminus[0,\infty)$. 
By Lemma \ref{inverse_of_rank_one_perturbation} and Lemma \ref{B_lambda_converges_to_B} {\it (ii)}, $A_\sigma$ is bijective with inverse 
\[
B_\sigma = B - \frac{\sigma}{\alpha_\sigma}\,\langle\,\cdot\,,\,B^\times p^1_K\,\rangle\,B1\nonumber
= B + \frac{\sigma}{R(R+\sigma)}\langle\,\cdot\,,\,VW_K\,\rangle W_K
\]
after simplifying using Lemma \ref{B_lambda_converges_to_B} {\it (iii)} and {\it (v)}. Now use
\[
\frac{\sigma}{R(R+\sigma)} = \frac{1}{R} - \frac{1}{R+\sigma}
                           = \frac{1}{R} + \sum_{k=-\infty}^{-1}\theta_k\eta^{k}
\]
with $\theta_k$ as in (\ref{formula_for_theta_k}). The expression (\ref{definition_of_B}) for $B$ leads to the result. 

\medskip

\noindent Now assume that $R=0$. In this case, $B_\sigma$ in (\ref{formula_for_theta_k}) becomes
\[
B_\sigma = B_0 -\frac{1}{\sigma}\langle\cdot, VW_K\rangle W_K.
\] 
For $u\in\mathscr{H}_{-s}$,
\begin{eqnarray*}
B_\sigma A_\sigma u & = & B_0 A u +\sigma\langle u, p^1_K\rangle B_0 1
-\frac{1}{\sigma}\langle Au, VW_K\rangle W_K - \langle u, p^1_K\rangle\langle 1, VW_K\rangle W_K\\
 & = & B_0Au + \langle u, p^1_K\rangle W_K
\end{eqnarray*}
by Lemma \ref{properties_of_B_0} {\it (i)} and Proposition \ref{properties_of_A_lambda} {\it (v)}. Each $u\in\mathscr{H}_{-s}$ may be written uniquely in the form
$
u = v + \alpha\,W_K
$
for some $v\in\mathscr{M}$ and $\alpha\in\mathbb{C}$. As $AW_K=0$ by Proposition \ref{properties_of_A_lambda} {\it (iv)}, we obtain
\[
B_\sigma A_\sigma u =  B_0Av + \alpha\,W_K = v + \alpha\,W_K = u
\]
by Lemma \ref{B_0_A_on_M} {\it (i)}. 
On the other hand,
\begin{eqnarray*}
A_\sigma B_\sigma u & = & A B_0 u +\sigma\langle B_0 u, p^1_K\rangle 1
-\frac{1}{\sigma}\langle u, VW_K\rangle AW_K - \langle u, VW_K\rangle\langle W_K,p^1_K\rangle 1\\
 & = & AB_0 u - \langle u, VW_K\rangle 1
\end{eqnarray*}
by Lemma \ref{properties_of_B_0} {\it (ii)} and Proposition \ref{properties_of_A_lambda} {\it (iv)}. Each $u\in\mathscr{H}_{-s}$ may be written uniquely in the form
$
u = w + \beta\,1
$
for some $w\in\mathscr{W}$ and $\beta\in\mathbb{C}$. So
\[
A_\sigma B_\sigma  u =  A B_0w - \beta\,\langle 1, VW_K\rangle 1 = 
w + \beta\,1= u
\]
by Lemma \ref{properties_of_B_0} {\it (i)} and Lemma \ref{B_0_A_on_M} {\it (ii)}. This shows that $B_\sigma$ is the inverse of $A_\sigma$ in the case $R=0$. 

\medskip

\noindent The final claim follows from the fact that
$
\sum_{k=-\infty}^{-1}\theta_k\eta^{k}=- \frac{1}{R+\sigma}
$
is bounded for small $\zeta$.
\qed

\begin{lemma}\label{technical_lemma_involving_B_sigma}
Let $l\in\mathbb{N}_0$ and $s>2l+1$. Then for small $\lambda>0$,
\begin{equation}\label{expansion_for_I_plus_B_V}
I+\sum_{j=1}^l\sum_{k=0}^1\lambda^j\eta^k B_\sigma A_j^k V = I - \sum_{j=1}^l\sum_{k=-\infty}^1\lambda^j\eta^{k}E_j^k
\end{equation}
in $B(\mathscr{H}_{-s})$. The coefficients $E_j^k$ are given by
\[
\begin{array}{llll}
E_j^{1} & = & - B_0A_j^1 V & \text{ for }j=1,2,\ldots,\\
E_j^0 & = & - B_0A_j^0 V - \theta_{-1}\,\langle\,\cdot\,,\,VA_j^{1}VW_K\,\rangle\,W_K & \text{ for }j=1,2,\ldots,\\
E_j^k & = & -\langle\,\cdot\,,V[\,\theta_k A_j^{0} + \theta_{k+1}A_j^{1}\,]VW_K\,\rangle\,W_K & \text{ for }j=1,2,\ldots\text{ and }k=-1,-2,\ldots.
\end{array}
\]
The double-summation in (\ref{expansion_for_I_plus_B_V}) converges absolutely in norm.
\end{lemma}

\noindent{\em Proof.}
Replace the expression for $B_\sigma$ as in Proposition \ref{expression_for_B_sigma} to obtain
\begin{eqnarray*}
I+\sum_{j=1}^l\sum_{k=0}^1\lambda^j\eta^k B_\sigma A_j^k V
& = &
I+\sum_{j=1}^l\sum_{k=0}^1\lambda^j\eta^k 
\left\{ 
B_0 + \sum_{p=-\infty}^{-1}\theta_p\eta^{p}\langle\,\cdot\,,\,VW_K\,\rangle\,W_K
\right\} A_j^k V
\end{eqnarray*}
\[
=
I+
\sum_{j=1}^l\sum_{k=0}^1\lambda^j\eta^k B_0A_j^k V
+
\sum_{j=1}^l\sum_{k=0}^1\sum_{p=-\infty}^{-1}\lambda^j\eta^{k+p}
\theta_p\langle\,\cdot\,,\,VA_j^{k}VW_K\,\rangle\,W_K
\]
\[
=
I+\sum_{j=1}^l\sum_{k=0}^1\lambda^j\eta^k B_0A_j^k V
+
\sum_{j=1}^l\sum_{p=-\infty}^{-1}\lambda^j\eta^{p}
\theta_p\langle\,\cdot\,,\,VA_j^{0}VW_K\,\rangle\,W_K
+
\sum_{j=1}^l\sum_{p=-\infty}^{-1}\lambda^j\eta^{1+p}
\theta_p\langle\,\cdot\,,\,VA_j^{1}VW_K\,\rangle\,W_K
\]
\begin{eqnarray*}
 & = & I+ \sum_{j=1}^l\sum_{k=0}^1\lambda^j\eta^k B_0A_j^k V\\
 &   & +\sum_{j=1}^l\lambda^j\theta_{-1}
 \langle\,\cdot\,,\,VA_j^{1}VW_K\,\rangle\,W_K\\
 &   & +\sum_{j=1}^l\sum_{k=-\infty}^{-1}\lambda^j\eta^{k}
 \left\{
 \theta_k\langle\,\cdot\,,\,VA_j^{0}VW_K\,\rangle\,W_K
+\theta_{k-1}\langle\,\cdot\,,\,VA_j^{1}VW_K\,\rangle\,W_K
 \right\}\\
  & = & I+ \sum_{j=1}^l\lambda^j\eta B_0A_j^1 V + \sum_{j=1}^l\lambda^j B_0A_j^0 V\\
 &   & +\sum_{j=1}^l\lambda^j\theta_{-1}
 \langle\,\cdot\,,\,VA_j^{1}VW_K\,\rangle\,W_K\\
 &   & +\sum_{j=1}^l\sum_{k=-\infty}^{-1}\lambda^j\eta^{k}
\langle\,\cdot\,,\,V\left[\, \theta_kA_j^{0}+
\theta_{k-1}A_j^{1}\,\right]VW_K\,\rangle\,W_K.
\end{eqnarray*}
From (\ref{formula_for_theta_k}) it can be seen that there exist finite constants $e_j$ ($j=1,2,\ldots$) such that
\begin{equation}\label{estimate_for_E_j_k}
\left|\left|\,E_j^k\,\right|\right|_{B(\mathscr{H}_{-s})}
\leq e_j\,\theta^{-k}
\end{equation}
for $k=\ldots,-2,-1$. In fact, an estimate of the above form also extends to the case $k=0$ and $k=1$. This shows that the double-summation converges absolutely in norm.
\qed

\begin{lemma}\label{technical_lemma_involving_inversion}
Let $l\in\mathbb{N}_0$ and $s>2l+1$. Then
\begin{equation}\label{expansion_for_I_minus_sum_inverse}
\left\{\,I - \sum_{j=1}^l\sum_{k=-\infty}^1
\lambda^j\eta^{k}E_j^k-o(\lambda^l)\,\right\}^{-1}
= I + \sum_{j=1}^l\sum_{k=-\infty}^j\lambda^j\eta^{k}D_j^k
+o(\lambda^l)
\end{equation}
in $B(\mathscr{H}_{-s})$ as $\lambda\downarrow 0$. The coefficients $D_j^k$ are given by
\[
D_j^k = \sum_{|\alpha|=j,\,|\beta|=k}E_\alpha^\beta
\]
where the multi-indices $(\alpha,\beta)$ belong to the set
\[
(\alpha,\,\beta)\in\bigcup_{n=1}^\infty\mathbb{N}^n\times\Lambda^n
\]
where $\Lambda:=\left\{\,\ldots,-2,-1,0,1\,\right\}$. The double-summation in (\ref{expansion_for_I_minus_sum_inverse}) converges absolutely in norm. 
\end{lemma}

\noindent{\em Proof.}
The operator
\[
T:=\sum_{j=1}^l\sum_{k=-\infty}^1\lambda^j\eta^{k}E_j^k+o(\lambda^l)
\]
satisfies $\left|\left|\,T\,\right|\right|_{B(\mathscr{H}_{-s})}<1$ for small $\lambda>0$. The inverse of $I-T$ may expressed as a Neumann series with coefficients as stated. 

\medskip

\noindent Suppose that $j\in\mathbb{N}$ and $\mathbb{Z}\ni k<0$. Using (\ref{estimate_for_E_j_k}),
\[
\left|\left|\,D_j^k\,\right|\right|_{B(\mathscr{H}_{-s})} \leq
\sum_{n=1}^j\sum_{|\alpha|=j}\sum_{|\beta|=k}\left|\left|\,E_\alpha^\beta\,\right|\right|_{B(\mathscr{H}_{-s})}
\leq
\theta^{-k}
\sum_{n=1}^j\sum_{|\alpha|=j}e_\alpha\,a(n,k).
\]
Assume that $k\leq -(3/2)j$. By Lemma \ref{lattice_point_counting_lemma}, the right-hand side may be estimated via
\[
\left\{\,\sum_{n=1}^j\sum_{|\alpha|=j}2^{n-1}\,a(n)\,e_\alpha\,\right\}
|k|^{j-1}\theta^{-k}.
\]
The index $n$ refers to the length of the multi-index $\alpha$. 
An inequality of the above form can be extended to the case $k<0$. In summary (for future use), for any $j\in\mathbb{N}$ and $k<0$,
\begin{equation}\label{estimate_for_D_j_k}
\left|\left|\,D_j^k\,\right|\right|_{B(\mathscr{H}_{-s})} \leq
d_j\,|k|^{j-1}\theta^{-k}
\end{equation}
for some finite constant $d_j$. This shows that the double-summation in (\ref{expansion_for_I_minus_sum_inverse}) converges absolutely in norm. 
\qed

\begin{lemma}\label{expansion_for_I_plus_R_V_inverse}
Let $l\in\mathbb{N}_0$ and $s>2l+1$. Then
\begin{equation}\label{formula_for_inverse_of_I_plus_R_V}
\left(\,I+R^{(-1)}(\mu-\imath\,0)V\,\right)^{-1}
= \sum_{j=0}^l\sum_{k=-\infty}^j\lambda^j\eta^{k}B_j^k
+o(\lambda^l)
\end{equation}
in $B(\mathscr{H}_{-s})$ as $\lambda\downarrow 0$.
The coefficients are given by
\begin{equation}\label{formulae_for_B_0_k}
\begin{array}{llll}
B_0^0 & = & B_0, & \\
B_0^k & = & \theta_k\,\langle\,\cdot\,,\,VW_K\,\rangle\,W_K & \text{ for }k=-1,-2,\ldots,\\
B_j^{j} & = & D_j^{j}B_0 & \text{ for }j=1,2,\ldots,\\
B_j^k & = & D_j^k B_0 +\sum_{p+q=k}\theta_p\,\langle\,\cdot\,,\,VW_K\,\rangle\,D_j^qW_K & \text{ for } j=1,2,\ldots\text{ and }k<j.
\end{array}
\end{equation}
The double-summation in (\ref{formula_for_inverse_of_I_plus_R_V}) converges absolutely in norm for small $\lambda>0$. 
\end{lemma}

\noindent{\em Proof.}
We rewrite (\ref{expansion_for_I_plus_R_V}) using Lemma \ref{technical_lemma_involving_B_sigma} as 
\begin{eqnarray*}
I+R^{(-1)}(\mu-\imath\,0)V
& = &
A_\sigma\,\left\{
I - \sum_{j=1}^l\sum_{k=-\infty}^1
\lambda^j\eta^{k}E_j^k-o(\lambda^l)
\right\}.
\end{eqnarray*}
Inverting using Lemmas \ref{technical_lemma_involving_inversion} and \ref{expression_for_B_sigma} we obtain
\begin{eqnarray*}
\left(\,I+R^{(-1)}(\mu-\imath\,0)V\,\right)^{-1} & = &
B_\sigma + \sum_{j=1}^l
\sum_{k=-\infty}^j
\lambda^j\eta^{k}D_j^kB_\sigma + o(\lambda^l)\\
 & = &
 B_0 + \sum_{k=-\infty}^{-1}\theta_k\eta^{k}\langle\cdot,VW_K\rangle W_K\\
 & & + \sum_{j=1}^l\sum_{k=-\infty}^j\lambda^j\eta^{k}D_j^k B_0\\
 & & + \sum_{j=1}^l\sum_{k=-\infty}^j\sum_{r=-\infty}^{-1}\lambda^j\eta^{k+r}\theta_r
 \langle\cdot, VW_K\rangle D_j^kW_K+ o(\lambda^l)\\
 & = &
 B_0 + \sum_{k=-\infty}^{-1}\theta_k\eta^{k}\langle\cdot,VW_K\rangle W_K\\
 & & + \sum_{j=1}^l\lambda^j\eta^jD_j^{j} B_0\\
 & & + \sum_{j=1}^l\sum_{k=-\infty}^{j-1}\lambda^j\eta^{k}
 \left\{\,D_j^k B_0 + \sum_{p+q=k}\theta_p\langle\cdot,VW_K\rangle D_j^q W_K\,\right\}
 + o(\lambda^l).
\end{eqnarray*}
For absolute convergence of the double-summation, let us first consider the term 
\[
\sum_{p+q=k,\,q<0}\theta_p\langle\cdot,VW_K\rangle D_j^q W_K
\]
for $j\in\mathbb{N}$ and $k\leq -2$. From (\ref{estimate_for_D_j_k}) its norm may be estimated by
\[
c\,\sum_{p+q=k,\,q<0}|\theta_p|\left|\left|\,D_j^q\,\right|\right|_{B(\mathscr{H}_{-s})}
\leq
c^\prime\,\left\{\sum_{p+q=k,\,q<0}|q|^{j-1}\right\}\theta^{-k}
\leq c^\prime\,|k|^{j}\theta^{-k}
\]
Consequently, for any $j\in\mathbb{N}_0$ and $k<0$,
\begin{equation}\label{estimate_for_B_j_k}
\left|\left|\,B_j^k\,\right|\right|_{B(\mathscr{H}_{-s})} \leq
b_j\,|k|^{j}\theta^{-k}
\end{equation}
for some finite constant $b_j$. So the double-summation converges absolutely in norm.
\qed

\begin{proposition}\label{expansion_of_T_of_lambda}
Let $l\in\mathbb{N}_0$. There exist $T_j^k\in\mathfrak{S}_1(\mathfrak{h})$, $0\leq j\leq 2l$, $-\infty< 2k\leq j$ such that
\begin{equation}\label{formula_for_T_of_lambda}
T(\lambda)=\sum_{0\leq j\leq 2l}\sum_{-\infty<2k\leq j}(\imath\,\lambda^{1/2})^j\eta^{k} T_j^k
+ o(\lambda^l)
\end{equation}
in $\mathfrak{S}_1(\mathfrak{h})$ as $\lambda\downarrow 0$. The coefficients are given by
\begin{equation}\label{formula_for_T_j_k}
T_j^k = 2\pi\imath\sum_{p+2q+r=j}(-1)^{q+r}U_p^{(-1)}VB_q^kU_r^{(-1)*}
\end{equation}
for $0\leq j\leq 2l$ and $-\infty< 2k\leq j$. The double-summation in (\ref{formula_for_T_of_lambda}) converges absolutely in norm. Also,
\[
T_0^0 = 0.
\]
\end{proposition}

\noindent{\em Proof.}
The argument proceeds as in \cite{McG2000} Proposition 4.4. Choose $s>2l+1$. As in
(\ref{expansion_for_U_minus_1}),
\[
U^{(-1)}(\mu) = \sum_{p=0}^{2l}(\imath\,\lambda^{1/2})^p U_p^{(-1)} + o(\lambda^l)
\]
in $\mathfrak{S}_2(\mathscr{H}_s,\mathfrak{h})$ as $\lambda\downarrow 0$. By Lemma \ref{expansion_for_I_plus_R_V_inverse},
\[
\left(\,I+R^{(-1)}(\mu-\imath\,0)V\,\right)^{-1}
= \sum_{q=0}^l\sum_{k=-\infty}^q(-1)^q(\imath\lambda^{1/2})^{2q}\eta^{k}B_q^k
+o(\lambda^l)
\]
in $B(\mathscr{H}_{-s})$ as $\lambda\downarrow 0$. The expansion follows straightforwardly.

\medskip

\noindent Fix $j\in\mathbb{N}_0$ and $k<0$. By (\ref{estimate_for_B_j_k}),
\[
\left|\left|\,
T_j^k
\,\right|\right|_{\mathfrak{S}_1(\mathfrak{h})}
\leq
c\,\sum_{p+2q+r=j}
\left|\left|\,B_q^k\,\right|\right|_{B(\mathscr{H}_{-s})}
\leq
c\,\sum_{p+2q+r=j}b_q |k|^{q}\theta^{-k}
\leq
c^\prime\,\left\{\sum_{p+2q+r=j}1\right\}|k|^{[j/2]}\theta^{-k}
\]
Thus for each $j\in\mathbb{N}_0$, there exists a finite constant $t_j$ such that
\begin{equation}\label{estimate_for_T_j_k}
\left|\left|\,T_j^k\,\right|\right|_{\mathfrak{S}_1(\mathfrak{h})}
\leq t_j \langle k\rangle^{[j/2]}\theta^{-k}
\end{equation}
for $-\infty<2k\leq j$. This establishes the summability claim.

\medskip

\noindent From (\ref{formula_for_T_j_k}), (\ref{formula_for_U_modified_0}) and (\ref{kernel_of_u_j}),
\[
T_0^0=2\pi\imath\,U_0VB_0U_0^*=(\imath/4\pi)\langle VB_01,1\rangle\langle\cdot\,,1\rangle.
\]
By Lemma \ref{properties_of_B_0} {\it (i)}, $B_01=0$; hence $T_0^0=0$.
\qed

\medskip

\begin{theorem}\label{asymptotic_expansion_of_xi}
Let $l\in\mathbb{N}$. Then 
\[
\xi(\lambda)=\sum_{k=-l}^{-1}
\xi_0^k\,\eta^k 
+ O(\eta^{-(l+1)})
\]
as $\lambda\downarrow 0$ where the coefficients are given by
\begin{equation}\label{expression_for_xi_j_k}
\xi_0^k=\frac{1}{2\pi\imath}\sum_{|\alpha|=0,\,|\beta|=k}\frac{(-1)^p}{p}\mathrm{Tr}\left[\,T_\alpha^\beta\,\right].
\end{equation}
\end{theorem}

\medskip

\noindent In the above, $p$ signifies the length of the multi-index $\alpha$ (resp. $\beta$). 

\noindent{\em Proof.}
For small $\lambda>0$,
\[
\xi(\lambda) = \frac{-1}{2\pi\imath}\,\mathrm{Tr}\log(I+T(\lambda))
\]
in virtue of (\ref{Birman_Krein_formula}). From Proposition \ref{expansion_of_T_of_lambda} we extract the expansion
\[
T(\lambda)=\sum_{k=-l}^{-1}\eta^{k} T_0^k
+ O(\eta^{-(l+1)})
\]
and insert into the formula
\[
\log(I+T) = \sum_{p=1}^\infty
\frac{(-1)^{p+1}}{p}\,T^p
\]
valid for $T\in B(\mathfrak{h})$ with $\left|\left| T\right|\right|<1$. 
\qed

\section{First three coefficients in low-energy expansion of the scattering phase}

\begin{lemma}\label{preliminary_formulae_for_xi_j_k}
The following identities hold:
\begin{itemize}
\item[{\it (i)}]
\[
\xi^{-1}_0 = - \frac{1}{2\pi\imath}\mathrm{Tr}[\,T^{-1}_0\,],
\]
\item[{\it (ii)}]
\[
\xi_0^{-2} = \frac{1}{2\pi\imath}\left\{\,-\mathrm{Tr}[\,T^{-2}_0\,] + (1/2)\,\mathrm{Tr}[\,T^{-1}_0 T^{-1}_0 \,]\,\right\},
\]
\item[{\it (iii)}]
\[
\xi_0^{-3} = \frac{1}{2\pi\imath}
\left\{\,
-\mathrm{Tr}[\,T^{-3}_0\,] 
+ (1/2)\,\left(\mathrm{Tr}[\,T^{-2}_0 T^{-1}_0 \,]+\mathrm{Tr}[\,T^{-1}_0 T^{-2}_0 \,]\right)
- (1/3)\,\mathrm{Tr}[\,T^{-1}_0 T^{-1}_0 T^{-1}_0\,]
\,\right\}.
\]
\end{itemize}
\end{lemma}

\noindent{\em Proof.}
These expressions follow directly from (\ref{expression_for_xi_j_k}).
\qed

\begin{theorem}\label{formulae_for_xi_0_etc}
The following identities hold:
\begin{itemize}
\item[{\it (i)}]
$\xi_0^{-1}=1$; 
\item[{\it (ii)}]
$\xi_0^{-2}=C(K) - \log\,4 + 2\gamma$; 
\item[{\it (iii)}]
$\xi_0^{-3}=\left(\,C(K) - \log\,4 + 2\gamma\right)^2 - \frac{\pi^2}{3}$.
\end{itemize}
\end{theorem}

\noindent{\em Proof.}
First note that from (\ref{formulae_for_B_0_k}) the identity
$
\langle 1, VB_0^{k}1\rangle = \theta_{k}
$
holds for any $\mathbb{Z}\ni k<0$. Moreover, from (\ref{formula_for_T_j_k}),
\[
T^{k}_0 
= 2\pi\imath\,U_0 V B^{k}_0 U_0^*
\]
for any $k<0$. With the help of \cite{McG2010} Corollary 7.2 {\it (i)} (or straightforwardly from (\ref{kernel_of_u_j})),
\begin{equation}\label{formula_for_trace_of_T_k_0}
\begin{array}{lll}
\mathrm{Tr}[\,T^{k}_0\,] & = & 2\pi\imath\,\mathrm{Tr}[U_0 V B^{k}_0U_0^*]\\
                         & = & 2\pi\imath\,(1/4\pi)\langle\,VB^{k}_01,1\,\rangle\\
                         & = & (\imath / 2 )\,\theta_{k}.
\end{array}
\end{equation}
From (\ref{formula_for_theta_k}) we have
\[
\begin{array}{lll}
\theta_{-1} & = & -\frac{1}{b},\\
\theta_{-2} & = & \frac{1}{b}\left(\frac{R+a}{b}\right),\\
\theta_{-3} & = & -\frac{1}{b}\left(\frac{R+a}{b}\right)^2
\end{array}
\]
with
\[
a=(1/2\pi)\,\left(\,\log\,2-\gamma\,\right)+\imath / 4\text{ and }b=1/4\pi.
\]

\medskip

\noindent{\it (i)} From (\ref{formula_for_trace_of_T_k_0}),
\[
\mathrm{Tr}[\,T^{-1}_0\,] = (\imath / 2 )\,\theta_{-1}= -2\pi\imath.
\]
This and Lemma \ref{preliminary_formulae_for_xi_j_k} {\it (i)} gives the first item.

\medskip

\noindent{\it (ii)} Using the above identity once more,
\[
\mathrm{Tr}\,[\, T^{-2}_0 \,] =  (\imath / 2)\,\theta_{-2}.
\]
With the help of \cite{McG2010} Corollary 7.2 {\it (vii)},
\begin{eqnarray*}
\mathrm{Tr}[T^{-1}_0T^{-1}_0] & = & (2\pi\imath)^2\mathrm{Tr}[U_0VB^{-1}_0U_0^*U_0VB^{-1}_0U_0^*]\\
 & = & - (2\pi)^2(1/4)(2\pi)^{-4}(2\pi)^2\langle VB^{-1}_0 1, 1\rangle^2\\
 & = & -(1/4)\,\theta_{-1}^2.
\end{eqnarray*}
By Lemma \ref{preliminary_formulae_for_xi_j_k} {\it (ii)},
\[
\xi_0^{-2} = \frac{1}{2\pi\imath}\left\{\,-(1/2)\,\imath\,\theta_{-2}- (1/8)\,\theta_{-1}^2\,\right\} 
=-4\pi\,\left[\, R + \Re\, a\,\right]
= C(K) - \log\,4 + 2\gamma.
\]

\medskip

\noindent{\it (iii)} From (\ref{formula_for_trace_of_T_k_0}),
\[
\mathrm{Tr}\,[\,T^{-3}_0\,] = (\imath / 2 )\,\theta_{-3}.
\]
By \cite{McG2010} Corollary 7.2 {\it (vii)},
\[
\begin{array}{lll}
\mathrm{Tr}[T^{-2}_0 T^{-1}_0] & = & (2\pi\imath)^2\,
\mathrm{Tr}[U_0 V B^{-2}_0U_0^* U_0 V B^{-1}_0U_0^*]\\
 & = & -(2\pi)^2(1/4)(2\pi)^{-4}(2\pi)^2
\langle VB^{-2}_0 1, 1\rangle\,
\langle VB^{-1}_0 1, 1\rangle\\
 & = & -(1 / 4 )\, \theta_{-2}\theta_{-1}.
\end{array}
\]
By \cite{McG2010} Corollary 7.2 {\it (viii)},
\[
\begin{array}{lll}
\mathrm{Tr}[T^{-1}_0 T^{-1}_0 T^{-1}_0] & = & 
(2\pi\imath)^3\,\mathrm{Tr}[U_0 V B^{-1}_0U_0^* U_0 V B^{-1}_0U_0^* U_0 V B^{-1}_0U_0^*]\\
 & = & -\imath(2\pi)^3(1/8)(2\pi)^{-6}(2\pi)^3\langle VB^{-1}_0 1, 1\rangle^3\\
 & = & - ( \imath / 8 )\, \theta_{-1}^3.
\end{array}
\]
By Lemma \ref{preliminary_formulae_for_xi_j_k} {\it (iii)} and some computation,
\[
\begin{array}{lll}
\xi_0^{-3} & = & \frac{1}{2\pi\imath}
\left\{\,
-(\imath/2)\theta_{-3}- (1/4)\theta_{-2}\theta_{-1} - (1/3)(-\imath/8)\theta_{-1}^3\,\right\}\\
 & = & (4\pi)^2\left\{\,(R+a-\imath/4)^2 - \frac{1}{48}\,\right\}\\
 & = & \left(\,C(K) - \log\,4 + 2\gamma\right)^2 - \frac{\pi^2}{3}.
\end{array}
\] 
\qed
\section{Asymptotics of the pinned Wiener sausage}

\medskip

\noindent We first remark that $\gamma(t)$ may be written purely analytically as
\begin{equation}\label{gamma_as_trace}
\gamma(t)=(4\pi t)\,\mathrm{Tr}\,[\,e^{-tH}-e^{-tH_Y}\,].
\end{equation}
Let $0<\delta<1$. For $k\in\mathbb{Z}$,
\begin{equation}\label{asymptotic_expansion_for_LT_of_powers_and_logs}
\int_0^\delta t e^{-t\lambda}(-\log\,\lambda)^k\,d\lambda
\sim
\sum_{r=0}^\infty
(-1)^r\binom{k}{r}\Gamma^{(r)}(1)(\log\,t)^{k-r}
\end{equation}
as $t\rightarrow\infty$ according to \cite{Erdelyi1960} Lemma 3. Recall that for $k<0$, the binomial is specified by 
\[
\binom{k}{r}=(-1)^r\binom{-k+r-1}{r}.
\]

\begin{theorem}
Let $l\in\mathbb{N}$. Then
\[
\gamma(t)=\sum_{k=-l}^{-1}\gamma_0^k\,t\,(\log\,t)^k + o(t\,(\log\,t)^{-l})
\]
as $t\rightarrow\infty$ where
\begin{equation}\label{formula_for_gamma_0_k}
\gamma_0^k=4\pi\sum_{s-r=k}\xi_0^s(-1)^r\binom{s}{r}\Gamma^{(r)}(1).
\end{equation}
The extra constraints $-\infty<s\leq -1$ and $r\geq 0$ apply in the summation. 
\end{theorem}

\noindent{\em Proof.}
We write
\[
\gamma(t) = (4\pi\,t)\left\{\,
\int_0^\delta te^{-t\lambda}\xi(\lambda)\,d\lambda
+\int_\delta^\infty te^{-t\lambda}\xi(\lambda)\,d\lambda
\right\}.
\]
In virtue of (\ref{integrability_of_xi}) the second term decays exponentially. Write
\[
\xi(\lambda)=\sum_{k=-l}^{-1}
\xi_0^k\,\eta^k 
+ O(\eta^{-(l+1)})
\]
according to Theorem \ref{asymptotic_expansion_of_xi}. By (\ref{asymptotic_expansion_for_LT_of_powers_and_logs}) the term
\[
(4\pi\,t)
\int_0^\delta te^{-t\lambda}\eta^{-(l+1)}\,d\lambda=o(t\,(\log\,t)^{-l})
\]
can be absorbed into the remainder. Again by (\ref{asymptotic_expansion_for_LT_of_powers_and_logs}), for $-l\leq k\leq-1$,
\[
(4\pi\,t)\int_0^\delta t e^{-t\lambda}\eta^k\,d\lambda
=
4\pi\sum_{r=0}^{k+l}(-1)^r\binom{k}{r}\Gamma^{(r)}(1)t(\log\,t)^{k-r}
+ o(t(\log\,t)^{-l})
\]
as $t\rightarrow\infty$. 
Therefore,
\begin{eqnarray*}
(4\pi\,t)\int_0^\delta te^{-t\lambda}\sum_{k=-l}^{-1}
\xi_0^k\,\eta^k\,d\lambda
& = &
\sum_{k=-l}^{-1}
\xi_0^k\,(4\pi\,t)\,\int_0^\delta te^{-t\lambda}\eta^k\,d\lambda\\
 & = & 
\sum_{k=-l}^{-1}\xi_0^k\,4\pi\sum_{r=0}^{k+l}(-1)^r\binom{k}{r}\Gamma^{(r)}(1)t(\log\,t)^{k-r}+ o(t(\log\,t)^{-l})\\
 & = & 
\sum_{k=-l}^{-1}\left\{\,4\pi\sum_{s-r=k}(-1)^r\xi_0^s\,\binom{s}{r}\Gamma^{(r)}(1)\,\right\}t(\log\,t)^{k}+ o(t(\log\,t)^{-l}).
\end{eqnarray*}
It is understood that $-\infty<s\leq -1$ and $r\geq 0$ in the summation.
\qed

\begin{corollary}
The following identities hold:
\begin{itemize}
\item[{\it (i)}]
$\gamma_0^{-1}=4\pi$;
\item[{\it (ii)}]
$\gamma_0^{-2}=4\pi\left\{\,C(K) + \gamma  - \log\,4 \right\}$;
\item[{\it (iii)}]
$\gamma_0^{-3}=4\pi\left\{\,(\,C(K)+\gamma - \log\,4)^2-\frac{\pi^2}{6}\right\}$.
\end{itemize} 
\end{corollary}

\noindent{\em Proof.}
From (\ref{formula_for_gamma_0_k}) we derive
\begin{itemize}
\item[{\it (a)}]
$\gamma_0^{-1}=4\pi\xi_0^{-1}$;
\item[{\it (b)}]
$\gamma_0^{-2}=4\pi\left\{\xi_0^{-1}\Gamma^{(1)}(1)+\xi_0^{-2}\right\}$;
\item[{\it (c)}]
$\gamma_0^{-3}=4\pi\left\{\xi_0^{-1}\Gamma^{(2)}(1)+2\xi_0^{-2}\Gamma^{(1)}(1)+\xi_0^{-3}
\right\}$.
\end{itemize} 
According to \cite{AbramowitzStegun1964} 6.4.2,
\[
\Gamma^{(1)}(1)=-\gamma,\hspace{1cm}\Gamma^{(2)}(1)=\gamma^2+\frac{\pi^2}{6}.
\]
The identities {\it (i)}-{\it (iii)} now follow with the help of Theorem \ref{formulae_for_xi_0_etc}.
\qed

\medskip

\noindent{\em Acknowledgement.}
I would like to express my thanks to Professor Michiel van den Berg for several helpful conversations on the topic of this paper and also for suggesting this problem in the first place. For this, I am also grateful to the anonymous referee of \cite{McG2010}.
\section{Appendix}

\noindent In this Appendix we prove Theorem \ref{expansion_of_R_of_zeta}, Lemma \ref{expansion_of_U_lambda} and Theorem \ref{extension_and_compactness_of_V}. 

\medskip

\begin{lemma}\label{boundedness_of_k}
The operator $K$ with convolution kernel
\[
k(x):=(\log\,|x|)\,|x|^\alpha
\hspace{5mm}
(-2<\alpha<\infty)
\]
belongs to $\mathfrak{S}_2(\mathscr{H}_s,\,\mathscr{H}_{-s})$ whenever $s>\alpha\vee 0 + 1$. 
\end{lemma}

\medskip

\noindent{\em Proof.}
Consider the operator $K$ with convolution kernel $k(x):=|x|^\alpha$. Suppose that $\alpha\geq 0$. For $s>\alpha+1$, 
\begin{eqnarray*}
\|K\|_{\mathfrak{S}_2(\mathscr{H}_s,\,\mathscr{H}_{-s})}^2 & = &
\int_{\mathbb{R}^2\times\mathbb{R}^2}\langle x\rangle^{-2s}|x-y|^{2\alpha}\langle y\rangle^{-2s}\,dy\,dx\\
& \leq & 4^{\alpha}\int_{\mathbb{R}^2\times\mathbb{R}^2}\langle x\rangle^{-2s+2\alpha}\langle y\rangle^{-2s+2\alpha}\,dy\,dx<\infty.
\end{eqnarray*} 
In case the kernel $k$ includes the logarithmic term, split the integral into a sum of integrals over the domains
\[
A_1:=\left\{(x,\,y)\in\mathbb{R}^2\times\mathbb{R}^2:\,0<|x-y|<1\right\}
\text{ and }
A_2:=\left\{(x,\,y)\in\mathbb{R}^2\times\mathbb{R}^2:\,|x-y|>1\right\}.
\]
On $A_1$ use Young's inequality (\cite{Davies1989} 1.1.4) and on $A_2$ use the inequality
\begin{equation}\label{estimate_on_logarithm}
\log\,|x-y|\leq 2^\varepsilon\varepsilon^{-1}\langle x\rangle^\varepsilon\langle y\rangle^\varepsilon
\end{equation}
valid for any $\varepsilon>0$. The above decomposition can also be used to deal with the case $-2<\alpha<0$. 
\qed

\medskip

\begin{lemma}\label{weight_function_over_complement_of_ball}
Let $\beta>2$. For $y\in\mathbb{R}^2$ and $0<r\leq 1/2$ define
\[
f(y,\,r):=\int_{|x-y|\geq r^{-1}}\langle x\rangle^{-\beta}\,dx.
\]
Then there exists a finite constant $c$ such that
\[
f(y,\,r)\leq
\left\{
\begin{array}{lcl}
c r^{\beta-2} & \text{ for } & |y|\leq\frac{1}{2r},\\
c             & \text{ for } & |y|>\frac{1}{2r}.
\end{array}
\right.
\]
\end{lemma}

\medskip

\noindent{\em Proof.}
The result for $|y|>1/2r$ is clear. Suppose that $|y|\leq 1/2r$. Then $1-r|y|\geq 1/2\geq r$. Thus, $B(0,\,r^{-1}-|y|)\subseteq B(y,\,r^{-1})$ and $B(0,\,1/2)\subseteq B(0,\,1-r|y|)$. This means that
\[
f(y,\,r)\leq\int_{|x|\geq r^{-1}-|y|}\langle x\rangle^{-\beta}\,dx
\leq r^{\beta-2}\int_{|x|\geq 1/2}|x|^{-\beta}\,dx.
\]
\qed

\medskip

\noindent Let $\varphi_1$ denote the indicator function of the interval $[0,\,1]$ and $\varphi_2:=1-\varphi_1$. 

\begin{lemma}\label{estimate_for_phi_2_K_zeta}
Let $K(\zeta)$ be the operator with convolution kernel
\[
k(x;\,\zeta):=\varphi_2(|\zeta|^{1/2}|x|)(\log\,|x|)\,|x|^\alpha
\hspace{5mm}
(\zeta\in\mathbb{C}\setminus[0,\,\infty),\,\alpha\in\mathbb{R}).
\]
Let $s>\alpha\vee 0 + 1$. Then $K(\zeta)$ belongs to $\mathfrak{S}_2(\mathscr{H}_s,\,\mathscr{H}_{-s})$ and
\[
\|K(\zeta)\|_{\mathfrak{S}_2(\mathscr{H}_s,\,\mathscr{H}_{-s})}
=O(|\zeta|^{(s-\alpha-1)/2})
\]
as $\zeta\rightarrow 0$. 
\end{lemma}

\medskip

\noindent{\em Proof.}
Consider the operator $K(\zeta)$ with convolution kernel $k(x;\zeta):=\varphi_2(|\zeta|^{1/2}|x|)|x|^\alpha$. Suppose that $\alpha\geq 0$. For $s>\alpha+1$ we find
\[
\|K(\zeta)\|_{\mathfrak{S}_2(\mathscr{H}_s,\,\mathscr{H}_{-s})}^2
\leq 4^{\alpha}\int_{\mathbb{R}^2}\langle y\rangle^{-\beta} f(y,\,|\zeta|^{1/2})\,dy
\]
where $f$ is defined as in Lemma \ref{weight_function_over_complement_of_ball} and $\beta:=2(s-\alpha)$ and $r:=|\zeta|^{1/2}$. Using the estimate in Lemma \ref{weight_function_over_complement_of_ball} this may be bounded by
\[
4^{\alpha}c\left\{
|\zeta|^{s-\alpha-1}\int_{|y|\leq 1/2|\zeta|^{1/2}}\langle y\rangle^{-\beta}\,dy
+
\int_{|y|>1/2|\zeta|^{1/2}}\langle y\rangle^{-\beta}\,dy
\right\}
\]
for $0<|\zeta|\leq 1/4$. The latter integral has order $O(|\zeta|^{s-\alpha-1})$ as $\zeta\rightarrow 0$. This gives the result for $\alpha\geq 0$. Now suppose that $\alpha<0$. For $s>1$,
\[
\|K(\zeta)\|_{\mathfrak{S}_2(\mathscr{H}_s,\,\mathscr{H}_{-s})}^2
\leq 
|\zeta|^{-\alpha}\int_{|x-y|\geq|\zeta|^{-1/2}}\langle x\rangle^{-2s}\langle y\rangle^{-2s}\,dy\,dx.
\]
Combining this with the result for $\alpha=0$ yields the result for this case. 
In case the kernel $k(\cdot;\zeta)$ includes the logarithmic term, make use of (\ref{estimate_on_logarithm}). 
\qed

\begin{lemma}\label{estimate_for_phi_1_K_zeta}
Let $K(\zeta)$ be the operator with convolution kernel
\[
k(x;\,\zeta):=\varphi_1(|\zeta|^{1/2}|x|)\,|x|^\alpha
\hspace{5mm}
(\zeta\in\mathbb{C}\setminus[0,\,\infty),\,\alpha>0).
\]
Let $\alpha<s\leq\alpha+1$ with $s>1$. Then $K(\zeta)$ belongs to $\mathfrak{S}_2(\mathscr{H}_s,\,\mathscr{H}_{-s})$ and
\[
\|K(\zeta)\|_{\mathfrak{S}_2(\mathscr{H}_s,\,\mathscr{H}_{-s})}
=
\left\{
\begin{array}{lcl}
O(|\zeta|^{(s-\alpha-1)/2})         & \text{ if } & \alpha<s<\alpha+1,\\
O((-\log\,|\zeta|)^{1/2})           & \text{ if } & s=\alpha+1,
\end{array}
\right.
\]
as $\zeta\rightarrow 0$.
\end{lemma}

\medskip

\noindent{\em Proof.}
Using the fact that
\[
|x-y|^{2\alpha}\leq 2^{2\alpha}\langle y\rangle^{2\alpha}\text{ for }|y|\geq |x|
\]
we have
\begin{eqnarray*}
\|K(\zeta)\|_{\mathfrak{S}_2(\mathscr{H}_s,\,\mathscr{H}_{-s})}^2 & = &
\int_{|x-y|\leq|\zeta|^{-1/2}}\langle x\rangle^{-2s}|x-y|^{2\alpha}\langle y\rangle^{-2s}\,dy\,dx\\
& \leq & 2^{2\alpha+1}\int_{\mathbb{R}^2}\langle x\rangle^{-2s}
\int_{|x-y|\leq|\zeta|^{-1/2}}\langle y\rangle^{2(\alpha-s)}\,dy\,dx\\
& \leq & 2^{2\alpha+1}\int_{\mathbb{R}^2}\langle x\rangle^{-2s}
\int_{|y|\leq |\zeta|^{-1/2}}\langle y\rangle^{2(\alpha-s)}\,dy\,dx.
\end{eqnarray*}
For $\alpha<s<\alpha+1$, 
\[
\int_{|y|\leq |\zeta|^{-1/2}}\langle y\rangle^{2(\alpha-s)}\,dy
\leq\frac{\pi}{\alpha-s+1}4^{\alpha-s+1}|\zeta|^{s-\alpha-1}
\]
for $0<|\zeta|<1$. On the other hand, for $s=\alpha+1$,
\[
\int_{|y|\leq |\zeta|^{-1/2}}\langle y\rangle^{2(\alpha-s)}\,dy
\leq 2\pi\left\{-\frac{1}{2}\log\,|\zeta| + \log 2\right\}
\]
again for $0<|\zeta|<1$. This leads to the result. 
\qed

\medskip

\noindent{\em Proof of Theorem \ref{expansion_of_R_of_zeta}}.
First recall that by \cite{AbramowitzStegun1964} 9.2.3,
\begin{equation}\label{estimate_for_H_1_0_of_z}
\left|\,H^{(1)}_0(z)\,\right|\leq c\,|\,z\,|^{-1/2}\text{ for }|z|>1\text{ and }
0<\mathrm{Arg}\,z<\pi.
\end{equation}
With $\varphi_1$, $\varphi_2$ as before set
\[
k^{(j)}(x;\,\zeta):=\varphi_j(|\zeta|^{1/2}|x|)k(x;\,\zeta)\hspace{1cm}(j=1,\,2)
\]
with $k(\cdot;\,\zeta)$ as in (\ref{expansion_for_resolvent_kernel_k}). Using Lemma \ref{boundedness_of_k} and the estimate (\ref{estimate_for_H_1_0_of_z}) for $k^{(2)}(x;\,\zeta)$ it may be seen that $R(\zeta)$ belongs to $\mathfrak{S}_2(\mathscr{H}_s,\,\mathscr{H}_{-s})$ for any $s>1$. 

\medskip

\noindent Let $l\in\mathbb{N}_0$ and $s>2l+1$. By Lemma \ref{boundedness_of_k} each of the operators $K^\varepsilon_j$ belongs to $\mathfrak{S}_2(\mathscr{H}_s,\,\mathscr{H}_{-s})$ for $j=0,\ldots,l$ and $\varepsilon=0,\,1$. Define 
\[
k_l(x\,;\zeta):=\sum_{j=0}^l\sum_{\varepsilon=0}^1\zeta^j\eta^\varepsilon\,k_j^\varepsilon(x)\hspace{1cm}(x\in\mathbb{R}^2\setminus\{0\})
\]
and the cut-off kernels $k_l^{(j)}(\cdot\,;\zeta)$ ($j=1,\,2$) as above. By (\ref{expansion_of_Hankel_function}) there exists a finite constant $c$ such that
\[
\big|k^{(1)}(x;\,\zeta) - k^{(1)}_l(x;\,\zeta)\big|
\leq c\,\varphi_1(|\zeta|^{1/2}|x|)|\zeta|^{l+1}\,|\eta|\,\Big(1+|\log\,|x|\,|\Big)|x|^{l+1}
\]
for small $\zeta$. Set $\alpha=l+1$. For $l\geq 1$ we have that $s>\alpha+1$. The remainder estimate follows from Lemma \ref{boundedness_of_k}. Consider the case $l=0$. If $s>2$ use Lemma \ref{boundedness_of_k}. If $1<s\leq 2$ use Lemma \ref{estimate_for_phi_1_K_zeta}. To deal with the logarithmic term consider the operators with kernels
\[
\varphi_j(|x|)
\left\{
k^{(1)}(x;\,\zeta) - k^{(1)}_l(x;\,\zeta)
\right\}
\hspace{1cm}(j=1,\,2).
\]
The operator corresponding to $j=1$ is bounded by Young's inequality. For the second use the fact that for any $\varepsilon>0$ there exists a finite constant $c$ such that
\[
\varphi_2(|x|)\log\,|x|\leq c\,|x|^{\varepsilon}\text{ for }x\in\mathbb{R}^2.
\]

\medskip

\noindent In view of (\ref{estimate_for_H_1_0_of_z}) we have that
\[
|k^{(2)}(x;\,\zeta)|\leq c\,|\zeta|^{-1/4}\varphi_2(|\zeta|^{1/2}|x|)|x|^{-1/2}
\]
and by Lemma \ref{estimate_for_phi_2_K_zeta} we obtain that \[
\|K^{(2)}(\zeta)\|_{\mathfrak{S}_2(\mathscr{H}_s,\,\mathscr{H}_{-s})}=O(|\zeta|^{\frac{s-1}{2}})=o(|\zeta|^l).
\]
Similar considerations can be used to deal with the terms in $K^{(2)}_l(\zeta)$. 
\qed

\medskip

\noindent{\em Proof of Lemma \ref{expansion_of_U_lambda}}.
The kernel of $U(\lambda)$ is given by
\[
u(\omega,\,x;\,\lambda)=\frac{1}{\sqrt{2}}(2\pi)^{-1}e^{-\imath\lambda^{1/2}\omega\cdot x}.
\]
Therefore, $u(\omega,\,x;\,\lambda)$ has an absolutely convergent series expansion of the form
\[
u(\omega,\,x;\,\lambda)=\sum_{j=0}^\infty(\imath\lambda^{1/2})^j u_j(\omega,\,x)
\]
where
\[
u_j(\omega,\,x):=\frac{1}{\sqrt{2}}(2\pi)^{-1}\frac{(-1)^j}{j!}(\omega\cdot x)^j.
\]
It is clear that $u(\omega,\,x;\,\lambda)$ is uniformly bounded. The truncated kernel will be written $u_l(\omega,\,x;\,\lambda)$. Let $\varphi_1$, $\varphi_2$ be as previously. Define
\[
u^{(j)}(\omega,\,x;\,\lambda):=\varphi_j(\lambda^{1/2}|x|)u(\omega,\,x;\,\lambda)
\hspace{1cm}
(j=1,2)
\]
and $u^{(j)}_l(\omega,\,x;\,\lambda)$ similarly. Denote the corresponding operators by $U^{(1)}(\lambda)$, etc. Let $U$ be the operator with kernel $u(\omega,\,x):=|x|^j$ ($j\in\mathbb{N}_0$). Note that $U\in\mathfrak{S}_2(\mathscr{H}_s,\,\mathfrak{h})$ if $s>j+1$. Let $U^{(j)}(\lambda)$ be the operator with kernel
\[
u^{(j)}(\omega, x;\lambda):=\varphi_j(\lambda^{1/2}|x|)|x|^{2j}.
\] 
Then $U^{(1)}(\lambda)$ has Hilbert-Schmidt norm $O(\lambda^{(s-j-1)/2})$ provided $s\leq j+1$. On the other hand, the operator $U^{(2)}(\lambda)$ has norm $O(\lambda^{(s-j-1)/2})$ if $s>j+1$.

\medskip

\noindent We have the estimate
\[
|u^{(1)}(\omega,x;\lambda) - u^{(1)}_l(\omega,x;\lambda)|
\leq c\,(\lambda^{1/2}|x|)^{l+1}.
\]
Thus $\|U^{(1)}(\lambda)-U^{(1)}_l(\lambda)\|_{
\mathfrak{S}_2(\mathscr{H}_s,\,\mathfrak{h})}=o(\lambda^{l/2})$ provided $s>l+1$. It is straightforward to see that the Hilbert-Schmidt norm of the difference $U^{(2)}(\lambda)-U^{(2)}_l(\lambda)$ admits an estimate of the same order in $\lambda$. 
\qed

\medskip

\noindent{\em Proof of Theorem \ref{extension_and_compactness_of_V}.}
The compactness statement is equivalent to the result $\langle\cdot\rangle^s V\langle\cdot\rangle^s\in\mathfrak{S}_\infty(\mathscr{H})$. Define
$V(t):=Je^{-tH_Y}J^* - e^{-tH}$. Then
\[
\langle\cdot\rangle^s V(2t)\langle\cdot\rangle^s 
=
\langle\cdot\rangle^s V(t)e^{-tH}\langle\cdot\rangle^s
+
\langle\cdot\rangle^s Je^{-tH_Y}J^*V(t)\langle\cdot\rangle^s.
\] 
The kernel $k(t;\,x,y)$ of $\langle\cdot\rangle^{-s} e^{-tH}\langle\cdot\rangle^s$ is well-known to be
\[
k(t;\,x,y)=(4\pi t)^{-1}\langle x\rangle^{-s} e^{-|x-y|^2/4t}\langle y\rangle^s.
\]
Using the inequality
\[
\langle y\rangle^s\leq 2^s(\langle x\rangle^{s}+\langle x-y\rangle^{s})
\]
we see that $k(t;\,x,y)$ is dominated by a square-integrable convolution kernel and hence by Young's inequality (\cite{Davies1989} 1.1.4 for example) $\langle\cdot\rangle^{-s} e^{-tH}\langle\cdot\rangle^s\in B(\mathscr{H},\,L^\infty(\mathbb{R}^2))$. Let $M=(\Omega,\mathscr{M},X_t,\mathbb{P}_{x})$ be Brownian motion on $\mathbb{R}^2$. Put $\sigma_K:=\inf\{t>0:\,X_t\in K\,\}$, the first hitting time of $K$. By the strong Markov property of Brownian motion,
\begin{eqnarray}\label{angle_V_t_angle_1}
\nonumber
|\langle\cdot\rangle^s V(t)\langle\cdot\rangle^s 1(x)| & = & \langle x\rangle^s
\mathbb{E}_x(\langle X_t\rangle^s:\,\sigma(K)<t)\\
 & = & \langle x\rangle^s\mathbb{E}_x(\mathbb{E}_{X_{\sigma(K)}}\langle X_{t-\sigma(K)}\rangle^s:\,\sigma(K)<t).
\end{eqnarray}
Since
\[
\sup_{y\in K}\sup_{0\leq\tau\leq t}\mathbb{E}_y\langle X_\tau\rangle^s
\]
is finite, (\ref{angle_V_t_angle_1}) is square-integrable. Therefore, $\langle\cdot\rangle^s V(t)\langle\cdot\rangle^s\in B(L^\infty(\mathbb{R}^2),\mathscr{H})$. We conclude that $\langle\cdot\rangle^s V(t)e^{-tH}\langle\cdot\rangle^s\in\mathfrak{S}_2(\mathscr{H})$ (see \cite{Stollmann1994} for example) and hence the same for $\langle\cdot\rangle^s V(2t)\langle\cdot\rangle^s$ by domination and duality. 

\medskip

\noindent Compactness of $\langle\cdot\rangle^s V\langle\cdot\rangle^s$ follows once we have shown that
\[
\int_0^\infty e^{-t}\|\langle\cdot\rangle^s V\langle\cdot\rangle^s\|_{B(\mathscr{H})}\,dt<\infty
\]
by \cite{Voigt1992} Theorem 1.3 and Remark 1.2 (b). Applying H\"older's inequality inside the functional integral we obtain for any $f\in\mathscr{H}$,
\[
\|\langle\cdot\rangle^s V(t)\langle\cdot\rangle^s f\|_2
\leq\sup_{y\in K}\sup_{0\leq\tau\leq t}(\mathbb{E}_y\langle X_\tau\rangle^{2s})^{1/2}
\sup_{x\in\mathbb{R}^2}\langle x\rangle^s\mathbb{P}_x(\sigma(K)<t)^{1/2}\|f\|_2.
\]
The known expression for the Brownian motion transition density yields that
\[
t\mapsto\sup_{y\in K}\sup_{0\leq\tau\leq t}(\mathbb{E}_y\langle X_\tau\rangle^{2s})^{1/2}
\] 
is $O(1)$ as $t\rightarrow 0+$ and $O(t^{s/2})$ as $t\rightarrow\infty$. The function
\[
t\mapsto\sup_{x\in\mathbb{R}^2}\langle x\rangle^s\mathbb{P}_x(\sigma(K)<t)^{1/2}
\]
has the same behaviour as can be seen using the "principle of not feeling the boundary". Thus the above integral is indeed finite. 
\qed

\end{document}